\RequirePackage{fix-cm}
\documentclass[smallextended]{svjour3}       
\smartqed  

\usepackage[latin1]{inputenc}
\usepackage{amsfonts}
\usepackage{graphicx}
\usepackage{epsfig}
\usepackage{amssymb,amsmath,color,theorem}
\usepackage{amsmath}
\usepackage{amssymb}
\usepackage{subfigure}

\newcommand\Estim[2]{\hat{D}^{#1}_{#2}}
\newcommand\Der[2]{D^{#1}_{#2}}
\newcommand\ErrEstim[2]{e^{#1}_{#2}}

\graphicspath{{E:/Multimedia/Images/Chercheurs/Dayan_Liu/Num_algorithms/}}

\begin{document}

\title{Error analysis of a class of derivative estimators for noisy signals}

\author{Da-yan Liu \and Olivier Gibaru \and Wilfrid Perruquetti}

\institute{D.Y. Liu  \at Laboratoire de Paul Painlevé,  Université  de Lille 1,
  59650, Villeneuve d'Ascq, France \\
              \email{dayan.liu@inria.fr}
           \and
           O. Gibaru \at Arts et Métiers ParisTech centre de
Lille, Laboratory of Applied Mathematics and Metrology (L2MA),  8
Boulevard Louis XIV, 59046 Lille Cedex, France
\\ \email{olivier.gibaru@ensam.eu}
  \and W. Perruquetti \at \'{E}cole Centrale de Lille, Laboratoire de LAGIS,
 BP 48, Cit\'e Scientifique, 59650
Villeneuve d'Ascq, France\\ \email{wilfrid.perruquetti@inria.fr}
\and D.Y. Liu \and O. Gibaru \and W. Perruquetti \at \'{E}quipe
Projet Non-A, INRIA Lille-Nord Europe
            Parc Scientifique
de la Haute Borne 40, avenue Halley B\^{a}t.A, Park Plaza, 59650
Villeneuve d'Ascq, France.}

\date{Received: date / Accepted: date}
\maketitle

\begin{abstract}
Recent algebraic parametric estimation
techniques (see \cite{garnier,mfhsr}) led to point-wise
derivative estimates by using only the iterated integral of
a noisy observation signal (see \cite{num0,num}).
In this paper, we extend such differentiation methods
by providing a larger  choice of parameters in these integrals: they can be reals. For this,  
the extension is done via a truncated Jacobi orthogonal series
expansion. Then, the noise error contribution of these derivative estimations is investigated: after proving the existence of such integral with a stochastic process noise, their statistical properties (mean value, variance and
covariance) are analyzed. In particular, the following important results are obtained:
\begin{description}
    \item[$a)$] the bias error term, due to the truncation, can be reduced by tuning the parameters,
    \item[$b)$] such estimators can cope with a large class of noises for which the mean and covariance are polynomials in time (with degree smaller than the order of derivative to be estimated),
  \item[$c)$] the variance of the noise error is shown to be smaller in the case of negative real parameters than it was in \cite{num0,num} for  integer values.
\end{description}
Consequently, these derivative estimations can be improved by tuning the parameters  according to the here obtained knowledge of the parameters' influence on the error bounds.


\end{abstract}

\section{Introduction}\label{sec_intro}

Numerical differentiation is concerned with the estimation of derivatives of
noisy time signals. This problem has attracted a lot of attention
from different points of view: observer design in the control
literature (see \cite{R6,R7,R18,R19,R22,R34}), digital filter in
signal processing (see \cite{R2,R5,R8,R29,R31}) and so on. 
The problem of numerical differentiation is ill-posed in the sense
that a small error in measurement data can induce a large error in
the approximate derivatives. Therefore, various numerical methods
have been developed to obtain stable algorithms more or less
sensitive to additive noise. They mainly fall into five classes: 
the finite difference methods \cite{I.R.Khan,R.Qu,A.G.Ramm}, the
mollification methods \cite{D.N.Hao,D.A.Murio,D.A.Murio2}, the
regularization methods \cite{G.Nakamura,T.Wei,Y.Wang}, the algebraic
methods \cite{num,num0} that are the roots of the here reported
results, the differentiation by integration methods
\cite{C.Lanczos,S.K.Rangarajana,Z.Wang,JCAM}, i.e. using the Lanczos
generalized derivatives.

The Lanczos generalized derivative $D_{T}x$, defined in
\cite{C.Lanczos} by
\begin{equation}\label{flanczos1}
\forall t_0 \in I, \ D_{T}x(t_0)
=\frac{3}{2T}\int
_{-1}^{1}\tau\,x(t_0+T\tau)\,d\tau,
\end{equation}
is an approximation to the first derivative of $x$ in the sense that
$D_{T}x(t_0)=x^{\prime}\left( t_0\right) + O(T^2) $, where $I$ is an
open interval of $\mathbb{R}$ and $2T>0$ is the length of the
integral window on which the estimates are calculated. It is aptly
called a method of \emph{differentiation by integration}.
Rangarajana and al. \cite{S.K.Rangarajana} generalized it  for
higher order derivatives with
\begin{equation}\label{flanczos2}
\forall t_0 \in I, \ D_{T}^{(n)}x(t_0)=\frac{1}{T^{n}}\int_{-1}^{1}\gamma_{n}L_{n}%
(\tau)\,x(t_0+T\tau)\,d\tau,\ n\in\mathbb{N},
\end{equation}
where $x$ is assumed to belong to $C^{n+2}(I)$  and $L_{n}$ is the
$n^{th}$ order Legendre polynomial defined on $[-1,1]$. The
coefficient $\gamma_{n}$ is equal to $\frac{1\times3\times5\times
\cdots\times(2n+1)}{2}$. By applying the scalar product of the
Taylor expansion of $x$ at $t_0$\ with $L_{n}$  they showed that
$D_{T}^{(n)}x(t_0)=x^{(n)}(t_0)+O(T^{2})$. By using
Richardson extrapolation Wang and al. \cite{Z.Wang} have improved
the convergence rate for obtaining high order Lanczos derivatives
with the following affine schemes for any $n\in\mathbb{N}$
\begin{equation}\label{flanczos3}
\forall t_0 \in I, \
D_{T,\lambda_{n}}^{(n)}x(t_0)=\frac{1}{T^{n}}\int_{-1}^{1}L_{n}(\tau)\left(
a_{n}\,x(t_0+T\tau)+b_{n}\,x(t_0+\lambda_{n}T\tau)\right)  d\tau,
\end{equation}
where $x$ is assumed to belong to $C^{n+4}(I)$, $a_{n}$, $b_{n}$ and
$\lambda_{n}$ are chosen such that $D_{T,\lambda_{n}}^{(n)}x(t_0)=x^{(n)}%
(t_0)+O(T^{4})$. Recently, Liu et al. \cite{JCAM} further reduced the convergence rate in these high order cases by using Jacobi polynomials in their derivative estimators. Let us mention that all these previous estimators
are given in the central case, $i.e.$ they use the interval
$[t_0-T,t_0+T]$ to estimate the derivative value $x^{(n)}(t_0)$.
Hence, these central estimators are only suited for off-line
applications whereas causal estimators using the interval
$[t_0-T,t_0]$ are well suited for on-line estimation which is of
importance in signal processing, automatic control and in general in many real time applications.

Very recently, Mboup, Fliess and Join  introduced in \cite{FLI04a_compression_cras} and analyzed in \cite{num,num0} a
new causal and anti-causal version of numerical differentiation by
integration method based on Jacobi polynomials
\begin{equation}\label{alien}
\forall t_0 \in I, \Estim{\mu,\kappa}{\beta
T} x^{(n)}(t_0)=\frac{\gamma^{\mu,\kappa}_{n}}{(\beta T)^{n}}\int_{0}^{1}
w^{\mu,\kappa}(\tau) P_n^{\mu,\kappa}(\tau) \,y(t_0+\beta
T\tau)\,d\tau,
\end{equation}
where $n\in\mathbb{N}$,  $y=x+\varpi$  is a noisy observation of $x$
which is assumed to be analytic, $\varpi$ denotes a noise,  $T >0$
is the length of the time window for integration, $\beta= -1$
(causal version) or $\beta= 1$ (anti-causal version) and
$P_n^{\mu,k}$ is the $n^{th}$ order Jacobi polynomial defined on
$[0,1]$ (see \cite{R}, \cite{R35}) by
\begin{equation}
\forall t \in [0,1],  \
P_{n}^{\mu,\kappa}(t)=\sum_{s=0}^{n}\binom{n+\mu}{s}\binom{n+\kappa}{n-s}\left(
t-1\right)  ^{n-s}\, t^{s},
\end{equation}
associated to the weight function
\begin{equation}\label{eq_weight}
w^{\mu,\kappa}(t)=t^{\kappa} (1-t)^{\mu}.
\end{equation}

Let us mention that originally, these estimators were
obtained  using an algebraic setting: for this the authors applied a differential operator on a truncation of the
Taylor series expansion of $x$ in the operational domain. This operator is given by
\begin{equation}
\begin{split}
\label{annihiator22} \Pi^{\mu,\kappa}_{n}=\frac{1}{s^{n+1+\mu}} \cdot
\frac{d^{n+\kappa}}{ds^{n+\kappa}}\cdot s^{n}, \ \text{ with } \kappa,\mu \in
\mathbb{N},
\end{split}
\end{equation}
$s$ being the Laplace variable. It is in fact an annihilator which ``kills'' the undesired terms except the one we want
to estimate (see \cite{num,num0} for more details). This is the reason why, originally in \cite{num,num0}, the parameters
\textbf{$\kappa, \mu$ were assumed to be integers ($\kappa, \mu,
\in\mathbb{N}$)}. In that case, the coefficients
\begin{equation}\label{eq_gamma}
\gamma^{\mu,\kappa}_{n}=\frac{n!(\mu+\kappa+2n+1)!}{(\mu+n)!(\kappa+n)!}.
\end{equation}
Let us emphasize that those methods,
which are algebraic and non-asymptotic, exhibit good robustness
properties with respect to corrupting noises, without the need of
knowing their statistical properties (see \cite{ans,shannon} for
more theoretical details). The robustness properties have already
been confirmed by numerous computer simulations and several
laboratory experiments.

The derivative estimators given by (\ref{alien}) contain two sources
of errors: the bias term error which comes from the truncation of
the Taylor series expansion  and the noise error contribution. Let us note that a precise analysis for the noise error contribution of a known noise has been done:
\begin{itemize}
\item in \cite{Med08}, for a specific identification method,
\item in \cite{Cdc09}, for discrete cases by using different integration methods,
\item in \cite{num}, which shows that an affine estimator induces a small time delay in the estimates while reducing the bias term error for \textbf{integer parameters $\kappa, \mu$}.
\end{itemize}
Thus, the aim of this paper is to \textbf{reduce the errors} of such derivative estimators.

\par For this, Section \ref{section02} allows $\kappa, \mu$ to be real:  since the use of (\ref{annihiator22}) induces some natural limitations ($\kappa, \mu,\in \mathbb{N}$) we use truncated Jacobi orthogonal series to obtain a natural extension of (\ref{alien}). Let us recall that, in \cite{num}, such truncated Jacobi orthogonal series were shown to be related to (\ref{annihiator22}) and  the obtained estimators (\ref{alien}). Such obtained estimators are thus called minimal Jacobi estimators and are clearly rooted from \cite{num}.

\par After providing such extension ($-1<\kappa, \mu \in \mathbb{R}$), Section 3 analyzes the bias term error: it is shown, for minimal
 Jacobi estimator, that if the $(n+1)^{th}$  derivative of the signal $x$ is slowly changing within the time window of observation then one can reduce the bias term error making  $T\frac{\kappa+n+1}{\mu+\kappa+2n+2}$ small by tuning the parameters $T,\kappa, \mu$. Lastly, it is shown that for affine estimator (real affine combination of the introduced minimal Jacobi estimators), the time delay  can be
reduced with the extended parameters comparing to the one obtained
in \cite{num}.

\par Section \ref{section03} analyzes the noise error contribution in order to provide some guides  for tuning the parameters $T,\kappa, \mu$. For this we consider two cases: continuous and discrete cases so as to give noise error bounds  by using the Bienaym\'{e}-Chebyshev inequality. In the first case for noises which are continuous parameter stochastic process with finite
second moments,  the mean value function and the covariance
kernel of these noise error contributions are calculated leading to:
\begin{itemize}
  \item for noises whose mean value function and covariance kernel are
polynomials of degree $r<n$ then the noise error contribution $\ErrEstim{\beta T}{\varpi}(t_0)=0  \text { almost surely}$,
  \item for Wiener or Poisson process some bounds are obtained for the noise error contribution: explicit for $n=1,2$ and it is shown how to deal with the general case.
\end{itemize}
The discrete case leads to similar results under some modified
assumptions.

\par In Section \ref{section04}, a parameter is introduced in order to reduce the
error due to a numerical integration method when the extended parameters become
negative. Then the comparisons of the Jacobi estimators with
extended parameters and the ones with original parameters are finally done for the cases where the noises are
respectively  a white Gaussian noise and a Wiener process noise.  The
integral of the total square error  and the classical $SNR$ are
considered for these comparisons.
Some interesting ``delay-free'' estimations in simulation
results are shown.

\section{Jacobi estimators}\label{section02}

\par  Let us  start with  $y=x+\varpi$  a noisy
observation on a finite time open interval $I \subset \mathbb{R}%
^{+}$ of a real valued smooth signal $x$ which $n^{th}$ derivative
has to be estimated ($n \in \mathbb{N}$), and $\varpi$ denotes a noise. Let  us assume that $x \in C^{n+1}(I)$. For any
$t_0 \in I$, we denote ${\cal D}_{t_0}=\{t\in \mathbb{R}^+; t_0+\beta t \in
I\}$ where $\beta=\pm 1$.

In the two following subsections, we aim at extending the parameters $\kappa,\mu$ used in the estimators described by (\ref{alien}) from $\mathbb{N}$ to $]-1,+\infty[$. To do so, we follow \cite{num} by taking the truncated Jacobi orthogonal series. Let us stress that (\ref{alien}) leads to two families of anti-causal ($\beta= 1$) and causal ($\beta= -1$) estimators.

\par  Let us mention that,  $\kappa\in \mathbb{R}$ means a non integer differentiation with respect to $s$ in (\ref{annihiator22}) but has nothing to do with the estimation of non integer derivatives of noisy time signals. However, such estimation of non integer derivatives could be tackled by similar technics.
\subsection{Minimal estimators}\label{sec_minimal_estimators}
Since $x \in C^{n+1}(I)$, we can take the Jacobi orthogonal series
expansion of $x^{(n)}$
\begin{equation}
\begin{split}
\label{x^n00}
\forall t_0 \in I,\  x^{(n)}(t_0+\beta T t)= \sum_{i \geq 0} \frac{\langle P_{i}^{\mu+n,\kappa+n}(\tau),x^{(n)}(\beta T\tau+t_0)%
\rangle}{\|P_{i}^{\mu+n,\kappa+n}\|^{2}} P_{i}^{\mu+n,\kappa+n}(t),
\end{split}
\end{equation}
where $t \in [0,1]$, $\beta=\pm 1$, $T \in {\cal D}_{t_0}$, $n \in
\mathbb{N}$, and $\kappa,\mu \in ]-1,+\infty[$.

By taking the first term in ($\ref{x^n00}$) with $t=0$, we get the
following estimations
\begin{equation}
\begin{split} \label{x^n_til}
\forall t_0 \in I, \ \Der{\mu,\kappa}{\beta T}x^{(n)}(t_0)&=\frac{\langle P_{0}^{\mu+n,\kappa+n}(\tau),x^{(n)}(\beta T\tau+t_0)%
\rangle}{\|P_{0}^{\mu+n,\kappa+n}\|^{2}} P_{0}^{\mu+n,\kappa+n}(0)\\
&=\frac{1}{B(\kappa+n+1,\mu+n+1)}\int_{0}^{1}
w^{\mu+n,\kappa+n}(\tau) \,x^{(n)}(t_0+\beta T\tau)\,d\tau,
\end{split}
\end{equation}
where $B(\cdot,\cdot)$ is the classical Beta function. Recall the
Rodrigues formula
\begin{equation}
\begin{split}
\label{derivee}
\frac{d^n}{d\tau^n}\left\{(1-\tau)^{\mu+n}\tau^{\kappa+n}\right\}=(-1)^n
n! (1-\tau )^{\mu}\tau^{\kappa} P_{n}^{\mu,\kappa}(\tau).
\end{split}
\end{equation}
Then, by taking $n$ times integration by parts and using the
Rodrigues formula in (\ref{x^n_til}) we get
\begin{equation}
\begin{split} \label{x^n_til2}
\forall t_0 \in I, \ \Der{\mu,\kappa}{\beta T}x^{(n)}(t_0)=\frac{\gamma^{\mu,\kappa}_{n}}{(\beta
T)^{n}}\int_{0}^{1}  w^{\mu,\kappa}(\tau) P_n^{\mu,\kappa}(\tau) \,x(t_0+\beta T\tau)\,d\tau,
\end{split}
\end{equation}
where $\gamma^{\mu,\kappa}_{n}=\frac{n!}{B(\kappa+n+1,\mu+n+1)}$
(the natural extension of (\ref{eq_gamma})). Now replacing $x$ in
(\ref{x^n_til2}) by its noisy observation $y$, one obtains
(\ref{alien})
with $\beta=\pm 1$, $T \in {\cal D}_{t_0}$, $n \in \mathbb{N}$, and
$\kappa,\mu \in ]-1,+\infty[$. Let us emphasis that these estimators were originally introduced in \cite{num} with
$\kappa, \mu \in \mathbb{N}$ and since they are obtained by taking the first term in the Jacobi series expansion, we call them \emph{minimal Jacobi estimators}. Hence, it is natural to extend these two parameters from $\mathbb{N}$ to $]-1,+\infty[$.

\subsection{Affine estimators}\label{sec_affine_estimators}
The Jacobi orthogonal series expansion of $x^{(n)}$
is the projection of $x^{(n)}$ onto the Jacobi orthogonal polynomial
basis. The minimal Jacobi estimators are introduced by taking the first term in the Jacobi series expansion of $x^{(n)}$ at $t=0$ in
the previous subsection. Let $q\in\mathbb{N}$ be the difference between $N$ (the truncation order of the Taylor series expansion of $x$) and $n$ (the order of derivative we want to estimate). We give in this subsection two families of estimators by taking the $q$  first terms
 in the Jacobi series expansion of $x^{(n)}$ at a set point $\xi$ as follows
\begin{equation}
\begin{split}
\label{x^n}
\forall t_0 \in I,\ \Der{\mu,\kappa}{\beta T,N,\xi}x^{(n)}(t_0): = \sum_{i = 0}^q \frac{\langle P_{i}^{\mu+n,\kappa+n}(\tau),x^{(n)}(\beta T\tau+t_0)%
\rangle}{\|P_{i}^{\mu+n,\kappa+n}\|^{2}} P_{i}^{\mu+n,\kappa+n}(\xi),
\end{split}
\end{equation}
where $\xi \in [0,1]$, $\beta=\pm 1$, $T \in {\cal D}_{t_0}$, $n \in
\mathbb{N}$, and $\kappa,\mu \in ]-1,+\infty[$.

These estimators were originally introduced in \cite{num} with  $\kappa,\mu \in \mathbb{N}$ and $q \leq n+\kappa$.
Moreover, it was shown that these estimators could be written as an affine combination of some minimal estimators.
Hence, (\ref{x^n}) proposes here two families of extended estimators which can also be written as follows

\begin{equation}
\begin{split} \label{16}
\Estim{\mu,\kappa}{\beta T,N,\xi}x^{(n)}(t_0)=
\displaystyle\sum_{l=0}^q {\lambda}_l(\xi)\,
\hat{D}^{\mu_l,\kappa_l}_{\beta T}x^{(n)}(t_0),
\end{split}
\end{equation}
where $(\kappa_l,\mu_l)=(\kappa+q+l,\mu+l) \in \mathbb{R}^2$ and the
minimal estimators $\hat{D}^{\mu_l,\kappa_l}_{\beta T}x^{(n)}(t_0)$
are defined by $(\ref{alien})$. The coefficients ${\lambda}_l(\xi)
\in \mathbb{R}$ are the same as the ones given
 in \cite{num}. Thus, we can call them \emph{affine Jacobi estimators}.

Let us denote by $p^{\mu,\kappa}_{\beta T, n}$  (respectively $p^{\mu,\kappa}_{\beta T, n, N, \xi}$) the power
functions used in the integral of the minimal estimators
$\Estim{\mu,\kappa}{\beta T}x^{(n)}(t_0)$ (respectively the affine estimators $\Estim{\mu,\kappa}{\beta T,N,\xi}x^{(n)}(t_0)$).
Then according to $(\ref{alien})$, we have
\begin{equation}
\begin{split}\label{minimal-power-function}
 p^{\mu,\kappa}_{\beta T, n}(\tau)&= \frac{\gamma^{\mu,\kappa}_{n}}{(\beta T)^{n}}
w^{\mu,\kappa}(\tau) P_n^{\mu,\kappa}(\tau).
\end{split}
\end{equation}
From  $(\ref{16})$, we can infer that
\begin{equation}
 p^{\mu,\kappa}_{\beta T, n, N, \xi}(\tau)=
\sum_{l=0}^q {\lambda}_l(\xi)\ p^{\mu_l,\kappa_l}_{\beta T, n}(\tau).\label{power-function}
\end{equation}

If we take $N=n$ and $\xi=0$ in the affine estimators, then we obtain
$\Estim{\mu,\kappa}{\beta T,N=n,\xi=0}x^{(n)}(t_0)=\Estim{\mu,\kappa}{\beta T}x^{(n)}(t_0)$.
Consequently $\Estim{\mu,\kappa}{\beta T,N,\xi}x^{(n)}(t_0)$
gives a general presentation for minimal estimators and affine
estimators. We call them \emph{Jacobi estimators}.

Now, by a direct adaptation from \cite{num},
we can obtain the following proposition by using some properties of the Jacobi orthogonal polynomials.

\begin{proposition}
Let $\Estim{\mu,\kappa}{\beta T}x^{(n)}(t_0)$ be the
minimal Jacobi estimators  with $n \geq 1$, then we have
\begin{equation}\label{estimateur0}
\Estim{\mu,\kappa}{\beta T}x^{(n)}(t_0)=-(A+B)\Estim{\mu,\kappa}{\beta T}x^{(n-1)}(t_0)+A\Estim{\mu,\kappa+1}{\beta
T}x^{(n-1)}(t_0)+B\Estim{\mu+1,\kappa}{\beta
T}x^{(n-1)}(t_0)
\end{equation}
where  $(n+\mu)A=-(n+\kappa)B=\frac{(\mu+\kappa+2n+1)(\mu+\kappa+2n)}{2
\beta T}$.
\end{proposition}

\noindent\textbf{Proof.} 
We need the following (adapted from
\cite{R} (p.782)):
\begin{equation}
\left(2n+2+{\mu}+{\kappa}\right)(1-\tau)\,P_{n}^{\mu+1,\kappa}(\tau)=
(1+n+\mu)P_{n}^{\mu,\kappa}(\tau)-(n+1)
P_{n+1}^{\mu,\kappa}(\tau).\label{19}
\end{equation}
\begin{equation}
\left(2n+2+{\mu}+{\kappa}\right)\tau \,P_{n}^{\mu,\kappa+1}(\tau)=
(1+n+\kappa)P_{n}^{\mu,\kappa}(\tau)+(n+1)
P_{n+1}^{\mu,\kappa}(\tau). \label{20}
\end{equation}
Subtracting (\ref{19}) from (\ref{20})  gives
\begin{equation}
\begin{split}\label{recurence}
P_{n+1}^{\mu,\kappa}(\tau)=\frac{\mu-\kappa}{2(n+1)}
P_{n}^{\mu,\kappa}(\tau)+\frac{2n+2+\kappa+\mu}{2(n+1)}\left[\tau
P_{n}^{\mu,\kappa+1}(\tau)-(1-\tau)P_{n}^{\mu+1,\kappa}(\tau)\right].
\end{split}
\end{equation}
By using (\ref{recurence}), $(\ref{alien})$ becomes
\begin{equation*}
\begin{split}
\Estim{\mu,\kappa}{\beta T}x^{(n)}(t_0)
&=\frac{\mu-\kappa}{2 n}\, c^{\mu,\kappa}_{\beta T,n}\int_{0}^{1}
(1-\tau )^{\mu}\tau^{\kappa} P_{n-1}^{\mu,\kappa}(\tau) \,y(\beta
T\tau+t_{0})d\tau\\
&+\frac{2 n+\kappa+\mu}{2 n}\, c^{\mu,\kappa}_{\beta T,n}\int_{0}^{1}
(1-\tau )^{\mu}\tau^{\kappa+1} P_{n-1}^{\mu,\kappa+1}(\tau)
\,y(\beta T\tau+t_{0})d\tau\\
&-\frac{2 n+\kappa+\mu}{2 n}\, c^{\mu,\kappa}_{\beta T,n}\int_{0}^{1}
(1-\tau )^{\mu+1}\tau^{\kappa} P_{n-1}^{\mu+1,\kappa}(\tau)
\,y(\beta T\tau+t_{0})d\tau.
\end{split}
\end{equation*}
with $c^{\mu,\kappa}_{\beta T,n}=\frac{\gamma^{\mu,\kappa}_{n}}{(\beta T)^{n}
}$.
Since
\begin{equation*}
\begin{split} c^{\mu,\kappa}_{\beta T,n}&=\frac{n!}{(\beta T)^{n}
}\frac{\Gamma(\mu+\kappa+2n+2)}{\Gamma(n+\kappa+1)\,\Gamma(\mu+n+1)}\\
&=\frac{n}{\beta T
}\frac{(\mu+\kappa+2n+1)(\mu+\kappa+2n)}{(n+\kappa)\,(\mu+n)}\, c^{\mu,\kappa}_{\beta T,n-1},\\
c^{\mu,\kappa}_{\beta T,n}&=\frac{n}{\beta T
}\frac{\mu+\kappa+2n+1}{n+\mu}\,c^{\mu,\kappa+1}_{\beta T,n-1},\\
c^{\mu,\kappa}_{\beta T,n}&=\frac{n}{\beta T
}\frac{\mu+\kappa+2n+1}{n+\kappa}\,c^{\mu+1,\kappa}_{\beta T,n-1}.
\end{split}
\end{equation*}
one can complete to obtain (\ref{estimateur0}).
 \hfill $\blacksquare$

\subsection{Two different sources of errors}
Since $\int_{0}^{1} p^{\mu,\kappa}_{\beta T, n}(\tau) d\tau=\frac{n!}{(\beta T)^n}$ and using the orthogonality of the Jacobi polynomials, one obtains
\begin{equation}
 \label{x_n}
{x}^{(n)}(t_0)=\int_{0}^{1}p^{\mu,\kappa}_{\beta T, n}(\tau) \,{x}_n(\beta T
\tau+t_0)\,d\tau, \forall t_0 \in I,
\end{equation}
where $p^{\mu,\kappa}_{\beta T, n}(\tau)$ is given by (\ref{minimal-power-function}) and ${x}_n(\beta T \tau+t_0)=\displaystyle\sum_{i=0}^n
\frac{\left(\beta T \tau\right)^i}{i !} x^{(i)}(t_0)$. Since
\begin{equation}\label{eq_rewrite_y}
x_n(\beta T\tau+t_{0})=y(\beta T\tau+t_{0})-R_n(\beta
T\tau+t_{0})-\varpi(\beta T\tau+t_{0})
\end{equation}
where  $R_n(\beta
T\tau+t_{0})= \frac{\left(\beta T \tau\right)^{n+1}}{(n+1)!}
x^{(n+1)}(\theta_\pm)$ with $\theta_- \in ]t_0-T\tau,t_0[$ and
$\theta_+ \in ]t_0,t_0+T\tau[$, we obtain by replacing (\ref{eq_rewrite_y}) in (\ref{x_n})
\begin{equation}
 \label{x_nbis}
{x}^{(n)}(t_0)=\Estim{\mu,\kappa}{\beta T}{x}^{(n)}(t_0)-(\ErrEstim{\mu,\kappa}{R_n,\beta
T}(t_0)+\ErrEstim{\mu,\kappa}{\varpi,\beta T}(t_0))
\end{equation}
where\begin{eqnarray}
       \ErrEstim{\mu,\kappa}{R_n,\beta
T}(t_0) &=& \int_{0}^{1}p^{\mu,\kappa}_{\beta T, n}(\tau) \,R_n(\beta T\tau+t_{0})\,d\tau, \\
       \ErrEstim{\mu,\kappa}{\varpi,\beta T}(t_0) &=& \int_{0}^{1}p^{\mu,\kappa}_{\beta T, n}(\tau) \,
\varpi(\beta T\tau+t_{0})\,d\tau.
     \end{eqnarray}
Thus the minimal
estimators $\Estim{\mu,\kappa}{\beta T}{x}^{(n)}(t_0)$ are
corrupted by two sources of errors:
\begin{itemize}
  \item  the bias
term errors $\ErrEstim{\mu,\kappa}{R_n,\beta
T}(t_0)$ which comes
from the truncation of the Taylor series expansion of   $x$,
  \item the noise error
contributions $\ErrEstim{\mu,\kappa}{\varpi,\beta T}(t_0)$.
\end{itemize}

For affine Jacobi estimator, as soon as $\displaystyle\sum_{l=0}^q
{\lambda}_l(\xi)=1$, using (\ref{16}) we have
\begin{eqnarray*}
  \Estim{\mu,\kappa}{\beta T,N,\xi}{x}^{(n)}(t_0) &=& {x}^{(n)}(t_0)+
\displaystyle\sum_{l=0}^q {\lambda}_l(\xi)\,
(\ErrEstim{\mu_l,\kappa_l}{R_n,\beta
T}(t_0)+\ErrEstim{\mu_l,\kappa_l}{\varpi,\beta T}(t_0)) \\
&=&{x}^{(n)}(t_0)+\ErrEstim{\mu,\kappa}{R_n,\beta T,N,\xi}(t_0)+\ErrEstim{\mu,\kappa}{\varpi,\beta
T,N,\xi}(t_0),
\end{eqnarray*}
where $\ErrEstim{\mu,\kappa}{R_n,\beta
T,N,\xi}(t_0)$ and $\ErrEstim{\mu,\kappa}{\varpi,\beta
T,N,\xi}(t_0)$ are respectively the bias term error and the noise
error contributions for these Jacobi estimators (minimal or not).
They are given by
\begin{equation}
\ErrEstim{\mu,\kappa}{R_n,\beta
T,N,\xi}(t_0)=\int_{0}^{1}p^{\mu,\kappa}_{\beta T, n,N,\xi}(\tau) \,
R_n(\beta T\tau+t_{0})d\tau, \label{bias error}
\end{equation}
\begin{equation}
\ErrEstim{\mu,\kappa}{\varpi,\beta
T,N,\xi}(t_0)=\int_{0}^{1}p^{\mu,\kappa}_{\beta T, n,N,\xi}(\tau) \
\varpi(\beta T\tau+t_{0}) d\tau. \label{noisy error}
\end{equation}


\section{Analysis of the bias error contribution}\label{sec_biais_error}
This analysis is done for minimal and affine Jacobi estimators. In both cases it is possible to reduce the bias term error by tuning the parameters
$T,\kappa,\mu$. For minimal Jacobi estimator one can get an overvaluation of this error using a Taylor expansion with integral reminder whereas in
 the affine case it is better to follow \cite{num}.

\subsection{Analysis for minimal Jacobi estimators}

The following result states that as soon the $n+1$ time derivative of the signal $x$ is slowly changing on the time window of observation then
one can reduce the bias term error making  $T\frac{\kappa+n+1}{\mu+\kappa+2n+2}$ as small as possible.

\begin{proposition} Let $\Estim{\mu,\kappa}{\beta T}{x}^{(n)}(t_0)$ be the minimal Jacobi estimators defined by (\ref{alien}) for
${x}^{(n)}(t_0)$. Then the corresponding bias terms errors can be
bounded by
\begin{equation} \label{bias error bound}
\begin{split}
C^{\mu,\kappa}_{T} I_{n+1}^+
 \leq &  \ErrEstim{\mu,\kappa}{R_{n},
T}(t_0) \leq C^{\mu,\kappa}_{T}  S_{n+1}^+, \\
C^{\mu,\kappa}_{-T} S_{n+1}^-
 \leq & \ErrEstim{\mu,\kappa}{R_{n},-
T}(t_0) \leq  C^{\mu,\kappa}_{-T} I_{n+1}^-,
\end{split}
\end{equation}
where $C^{\mu,\kappa}_{\beta T}=\beta T\frac{\kappa+n+1}{\mu+\kappa+2n+2}$ and
\begin{equation} \label{bias error cond}
\begin{split}
I_{n+1}^+=\inf_{t_0 <\theta_+ < t_0+T} {x}^{(n+1)}(\theta_+),
S_{n+1}^+=\sup_{t_0 <\theta_+ < t_0+T}
{x}^{(n+1)}(\theta_+),\\
S_{n+1}^-= \sup_{t_0-T <\theta_- < t_0} {x}^{(n+1)}(\theta_-),
I_{n+1}^-= \inf_{t_0-T <\theta_- < t_0} {x}^{(n+1)}(\theta_-).
\end{split}
\end{equation}
\end{proposition}

\noindent\textbf{Proof.} Let us take the Taylor series expansion of
$x^{(n)}(t_0+\beta T\tau)$ in (\ref{x^n_til}), then we have
\begin{equation*}
\begin{split}
\forall t_0 \in I, \ \Estim{\mu,\kappa}{\beta T}{x}^{(n)}(t_0)
&=\gamma^{\mu,\kappa}_{n}\int_{0}^{1} w^{\mu,\kappa}(\tau)
\,x^{(n)}(t_0+\beta T\tau)\,d\tau \\
&=\gamma^{\mu,\kappa}_{n}\int_{0}^{1} w^{\mu,\kappa}(\tau)
\left(x_n^{(n)}(t_0+\beta T\tau)+ \hat{R}_n(\beta
T\tau+t_{0})\right) \,d\tau,
\end{split}
\end{equation*}
 where  $\hat{R}_n(\beta
T\tau+t_{0})=\beta T \tau  x^{(n+1)}(\hat{\theta}_\pm)$ with
$\hat{\theta}_- \in ]t_0-T\tau,t_0[$ and $\hat{\theta}_+ \in
]t_0,t_0+T\tau[$. Thus, the bias term errors is given by
\begin{equation}
\begin{split}
\ErrEstim{\mu,\kappa}{R_{n},\beta T}(t_0) &=\beta T
\gamma^{\mu,\kappa}_{n}\int_{0}^{1} w^{\mu,\kappa+1}(\tau)
 x^{(n+1)}(\hat{\theta}_\pm) \,d\tau.
\end{split}
\end{equation}
Then, this proof can be easily completed by taking the Beta function
and the  extreme  values of $x^{(n+1)}(\hat{\theta}_\pm)$. \hfill
$\blacksquare$

As shown in \cite{num}, when $\beta=-1$ (causal case),
$-C^{\mu,\kappa}_{-T}=C^{\mu,\kappa}_{T}=\frac{\kappa+n+1}{\mu+\kappa+2n+2}T$ is the time delay
when we estimate ${x}^{(n)}(t_0)$ by
$\Estim{\mu,\kappa}{-T}{x}^{(n)}(t_0)$.

\begin{corollary}
If $I_{n+1}^+\simeq S_{n+1}^+$ and $I_{n+1}^-\simeq S_{n+1}^-$
then  by minimizing this
delay $C^{\mu,\kappa}_{T}$ we also minimize the bias term errors.
\end{corollary}

\begin{figure}[h!]
\begin{center}
\includegraphics[scale=0.6]{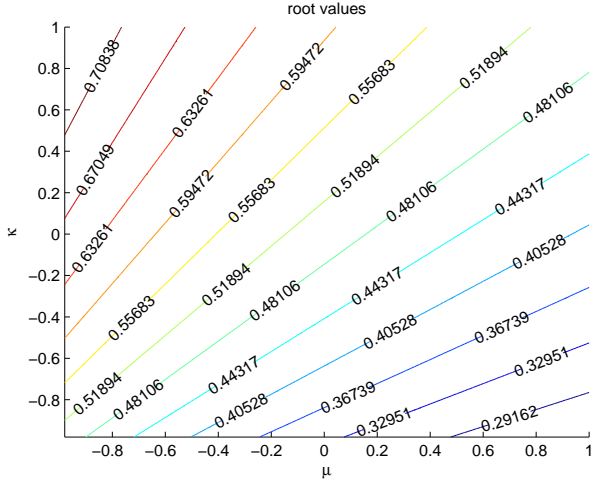}
 \caption{ Variation of $\frac{\kappa+2}{\kappa+\mu+4}$
with respect to $\kappa$ and $\mu$ .}
\label{root0}
\end{center}
\end{figure}

Since, when $\kappa,\mu \in ]-1,+\infty[$, $\frac{\kappa+n+1}{\mu+\kappa+2n+2}$ increases with respect to
$\kappa$ and decreases  with respect to $\mu$, the negative values
of $\kappa$ produce smaller bias term errors than the ones produced
by integer values of $\kappa$. This is one of the reason to extend the values of
$\kappa$.
It is clear that one can achieve a given bias term error by increasing $\mu$ and reducing $T$ (even choosing $\kappa,\mu$ as integer) but,
as we will see later on in Section \ref{section03}, it will increase the variance of the noise error contribution.
 When $n=1$, we can
see the variation of $\frac{\kappa+2}{\kappa+\mu+4}$  with respect to $(\kappa,\mu) \in ]-1,1]^2$ in
Figure $\ref{root0}$.

\subsection{Analysis for affine Jacobi estimators}

\begin{figure}[h!]
\centering {\includegraphics[scale=0.6]{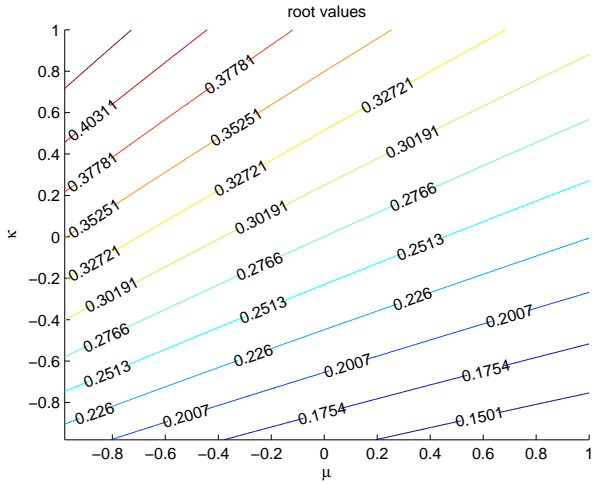}} \caption{ Variation
of $\xi(\kappa,\mu,n=1)$ with respect to $\kappa$ and $\mu$.}
\label{root}
\end{figure}

It was shown in \cite{num} (for $\beta=-1$ causal case), that the
bias term error $\ErrEstim{\mu,\kappa}{R_{n},-T,N,\xi}(t_0)$ in
(\ref{16})
produces a time delay of value $T\xi$, this is (when  there is no
noise $\varpi=0$):
$$\Estim{\mu,\kappa}{-T,N,\xi}{x}^{(n)}(t_0)\approx {x}^{(n)}(t_0-T\xi).$$
We always take the  value of $\xi$ as the
smallest root of the Jacobi polynomial $P_{q+1}^{\mu+n,\kappa+n}$
$(q=N-n)$, such that this affine causal estimator may be
significantly improved by admitting the minimal time delay. Hence,
$\xi$ is a function of $\kappa$, $\mu$ and $n$. We denote it by
$\xi(\kappa,\mu,n)$. We can see the variation of $\xi(\kappa,\mu, n=1)$ with respect to $(\kappa,\mu) \in ]-1,1]^2$ in Figure
$\ref{root}$. Hence, the extended parameters values give smaller
value for $\xi(\kappa,\mu,n=1)$.

\section{Analysis of the noise error contribution}\label{section03}

\par Before analyzing the resulting noisy error (\ref{noisy error}), let us study the existence of the integrals in the expressions of Jacobi
estimators. As the noisy observation $y$ is the sum of  $x$ and
noise $\varpi$, the Jacobi estimators are well defined if and only
if noise $\varpi$ is integrable. Indeed, according to
(\ref{minimal-power-function}) and (\ref{power-function}), if
$\varpi$ is an integrable function then  integrability of
$w^{\mu,\kappa}(\tau) \varpi (t_0+ T\tau)$ holds for $\mu,\kappa \in
]-1, +\infty[$ and $T \in {\cal D}_{t_0}$. Thus, in that case, the
integrals in the Jacobi estimators exist. Now, if $\varpi$ is a
continuous parameter stochastic process (see \cite{SP}) the next
result (Lemma \ref{stoc}) proves the existence of these integrals
and thus justifies (\ref{alien}), (\ref{16}) and (\ref{noisy error})
 as soon as the integrals are understood in the sense of convergence in
mean square (see Proposition \ref{corollary0}). For this, the
stochastic process $\{\varpi(\tau), \tau \geq 0\}$ should   satisfy
the following condition
\begin{description}
  \item[$(C_1):$] $\{\varpi(\tau), \tau \geq 0\}$ is   a continuous parameter stochastic process with finite second
moments, whose mean value function and covariance kernel are
continuous functions.
\end{description}

\begin{lemma} \label{stoc}
Let $\{\varpi(\tau), \tau \geq 0\}$ be a stochastic process
satisfying condition $(C_1)$. Then for any $t_0 \in I$ and $T \in {\cal D}_{t_0}$, the integral $\int_0^1
w^{\mu,\kappa}(\tau)\varpi (t_0+ T\tau) d\tau$ (with
$\mu,\kappa \in ]-1, +\infty[$) is well defined as a limit in
mean square of the usual approximating sum of the following form
\begin{equation}
\begin{split}
\label{riemann} \int_0^1 Y (\tau) d\tau=\lim_{m \rightarrow \infty}
\sum_{l=1}^m (\tau_l-\tau_{l-1})\,Y_l,
\end{split}
\end{equation}
where $Y(\tau)=w^{\mu,\kappa}(\tau) \varpi (t_0+
T\tau)$, $Y_l=Y(\xi_l)$ for any $\xi_l \in ]\tau_{l-1},\tau_l[$ and
$0=\tau_0<\tau_1<\cdots<\tau_m=1$ is a subdivision of the interval
$]0,1[$, such that $\displaystyle\max_{\substack{l=1,\cdots,m}}
(\tau_l-\tau_{l-1})$ tends to $0$ when $m$ tends to infinite.
\end{lemma}

\noindent\textbf{Proof of Lemma \ref{stoc}.}
For any fixed $t_0 \in D$, it was shown in \cite{loeve} (p. 472)
that if $\{Y(\tau), 0 <\tau <1\}$, where $Y(\tau)=w^{\mu,\kappa}(\tau) \varpi (t_0+T\tau)$, is a continuous parameter
stochastic process with finite second moments, then a necessary and
sufficient condition such that the family of approximating sums on
the right-hand side of $(\ref{riemann})$ has a limit in the sense of
convergence in mean square is
 that the double integral
$\int_0^1 \int_0^1 E[Y(s) Y(\tau)] \,d s \,d\tau$ exists.\\
Since for any $\tau \in ]0,1[$, $(1-\tau)^{\alpha} \, \tau^{\beta} <
\infty$, and $\{\varpi(\tau), \tau \geq 0\}$ is a continuous
parameter stochastic process with finite second moments, so does
$\{Y(\tau), 0<\tau<1\}$ for any $t_0 \in I$. Moreover, since the
mean value function and covariance kernel of $\varpi(\tau)$ are
continuous functions, so does $E[\varpi (t_0+ T\tau)\, \varpi (t_0+
Ts)]$ for all $\tau, s \in [0,1]$. Hence, $E[\varpi (t_0+ T\tau)\,
\varpi (t_0+ Ts)]$ is bounded for all $\tau, s \in [0,1]$.

Consequently, $\int_0^1 \int_0^1 \,
w^{\mu,\kappa}(\tau)w^{\mu,\kappa}(s)\,E[\varpi (t_0+ T\tau)\,
\varpi (t_0+ Ts)] \,d s \,d\tau $ exists when $\kappa,\mu \in ]-1,
+\infty[$, which implies that $(\ref{riemann})$ holds.\hfill
$\blacksquare$

If we take $y$  instead of $\varpi$ in the previous lemma, then we
can obtain the following proposition.

\begin{proposition} \label{corollary0} If
$x \in C^{n+1}(I)$, and the noise $\varpi$ satisfies condition $(C_1)$, then for any $t_0 \in I$, the integrals in the Jacobi estimators exist
in the sense of convergence in mean square.
\end{proposition}

From now on, we can investigate the noise error contributions for
these Jacobi estimators. Mainly the Bienaym\'{e}-Chebyshev
inequality is used to give two error bounds for these errors. Let us
denote the noise error contributions for the Jacobi estimators
$\ErrEstim{\mu,\kappa}{\varpi,\beta T,N,\xi}(t_0)$ (see (\ref{noisy
error})) by $\ErrEstim{\beta T}{\varpi}(t_0)$, then for any real
number $\gamma>0$
\begin{equation} \label{B-Chebyshev}
Pr\left(\left|\ErrEstim{\beta T}{\varpi}(t_0)-E[\ErrEstim{\beta T}{\varpi}(t_0)]\right| < \gamma\sqrt{Var[\ErrEstim{\beta T}{\varpi}(t_0)]}\right) > 1-
\frac{1}{\gamma^2},
\end{equation}
$i.e.$ {the probability for } $\ErrEstim{\beta T}{\varpi}(t_0)$ { to be
within the interval  } $\left] M_l \,, M_h\right[$ { is higher than
} $1- \frac{1}{\gamma^2}$, where $M_l=E[\ErrEstim{\beta T}{\varpi}(t_0)]-
\gamma \sqrt{Var[\ErrEstim{\beta T}{\varpi}(t_0)]}$ and
$M_h=E[\ErrEstim{\beta T}{\varpi}(t_0)]+\gamma
\sqrt{Var[\ErrEstim{\beta T}{\varpi}(t_0)]}$. These error bounds $M_l$, $M_h$ depend on the  parameters $\kappa$, $\mu$, $T$
and $\xi$ which can help us  in minimizing the noise error
contributions. From the previous section, we have extended the values of
$\kappa, \mu$ from $\mathbb{N}$ to $]-1,+\infty[$. Hence, we obtain
a higher degree of freedom so as to minimize the noise effects on
our estimators. In order to obtain these bounds we need to compute the means and variances of these
errors.

To do so, firstly, Subsection \ref{subsection_stochastic} considers
noises as continuous stochastic processes: it is shown that such
Jacobi estimators can cope with a large class of noises with mean
and covariance polynomials in time for which $M_l,M_h$ are obtained.
Let us note that this class includes well known processes such as
the Wiener  and the Poisson ones. Secondly, Subsection
\ref{subsection_discrete} deals with the discrete case for these
noises.

\subsection{Noise error contribution in the context of a stochastic process noise}\label{subsection_stochastic}

Let us assume that noise $\varpi$ satisfies condition $(C_1)$. To
simplify our notations, let us denote the power functions
$p^{\mu,\kappa}_{\beta T, n,N,\xi}$ associated to the Jacobi
estimators by $p^{\beta T}$. Then by applying   Theorem 3A in
\cite{SP} (p.79)  the means, variances  and  covariances of the
noise error contributions for the Jacobi estimators are given as
follows ($\forall \, T>0, T_1>0, \, {\rm and}\, T_2>0$)
\begin{equation}\label{mean2}
E\left[e_\varpi^{\beta T}(t_0)\right]=
 \int_0^1 p^{\beta T}(\tau)\ E\left[\varpi(t_0 +\beta T\tau)\right] \,d\tau,
\end{equation}
\begin{equation}\label{covariance2}
\begin{split}
&Cov\left[e_\varpi^{\beta T_1}(t_0),e_\varpi^{\beta T_2}(t_0)
\right]\\
&=\int_0^1 \int_0^1 p^{\beta T_1}(s)\,p^{\beta T_2}(\tau)\,
Cov\left[ \varpi(t_0+\beta T_1 s),\, \varpi(t_0+\beta T_2
\tau)\right] \,ds \,d\tau,\\
\end{split}
\end{equation}
\begin{equation}\label{var2}
Var\left[e_\varpi^{\beta T}(t_0)\right]=Cov\left[e_\varpi^{\beta
T}(t_0),e_\varpi^{\beta T}(t_0)\right].
\end{equation}

By using the property of  power function $p^{\mu,\kappa}_{\beta T,
n}$ defined in $(\ref{minimal-power-function})$, the following
theorem shows that such Jacobi estimators can deal with a large
class of noises for which the mean and covariance are polynomials in
time satisfying  the following conditions
\begin{description}
    \item[$(C_2):$] $ \forall (t_0+\tau) \in I$, the following holds
\begin{equation}
E[\varpi(t_0+\tau)]=  \sum_{i=0}^{n-1} \nu_i \, t_0^{k_1(i)}
\tau^i+E[\varpi(\tau)], \label{condition1}
\end{equation}
\begin{equation}
Cov[\varpi(t_0+s),\varpi(t_0+\tau)]=\displaystyle\sum_{i=0}^{n_1}
\eta_i \,t_0^{k_2(i)} \tau^i\,\displaystyle\sum_{i=0}^{n_2} \eta'_i \,t_0^{k_3(i)}
s^i+Cov[\varpi(s),\varpi(\tau)], \label{condition2}
\end{equation}
where ${k_1(i)} \in \mathbb{N}$, ${k_2(i)} \in \mathbb{N}$, ${k_3(i)} \in \mathbb{N}$, $\nu_i \in \mathbb{R}$, $\eta_i \in \mathbb{R}$, $\eta'_i \in
\mathbb{R}$ and $n_1 \in \mathbb{N}$, $n_2 \in \mathbb{N}$ such that
$\min(n_1,n_2) \leq n-1$.
\item[$(C_3):$] $\forall \tau \in I$, the following holds
\begin{equation}
E[\varpi(\tau)]=\displaystyle\sum_{i=0}^{n-1}\bar{\nu}_i \,\tau^i,
\end{equation}
\begin{equation}
Cov[\varpi(s),\varpi(\tau)]=\displaystyle\sum_{i=0}^{n_1} \bar{\eta}_i \,\tau^i\,\displaystyle\sum_{i=0}^{n_2} \bar{\eta}'_i \,s^i,
\end{equation}
where $\bar{\nu}_i \in \mathbb{R},\bar{\eta}_i, \bar{\eta}'_i \in \mathbb{R}$ and $\min(n_1,n_2)\leq n-1$
\end{description}

\begin{theorem} \label{corollary}
Let $\ErrEstim{\beta T}{\varpi}(t_0)$ be the noise error contribution for
the Jacobi estimator $\Estim{\mu,\kappa}{\beta T,N,\xi} {x}^{(n)}({t_0})$ where the noise  $\{\varpi(\tau), \tau \geq 0\}$ satisfies
 conditions $(C_1)$ and
 $(C_2)$. If $n \in \mathbb{N}^*$, then the mean, variance and covariance of
$e_\varpi^{\beta T}(t_0)$ do not depend on $t_0$. If in addition the noise $\{\varpi(\tau), \tau \geq 0\}$ satisfies conditions $(C_3)$
then $E[\ErrEstim{\beta T}{\varpi}(t_0)]=0$, $Cov[e_\varpi^{\beta T_1}(t_0),e_\varpi^{\beta T_2}(t_0)]=0$ and
$Var[\ErrEstim{\beta T}{\varpi}(t_0)]=0$.
\end{theorem}

\noindent\textbf{Proof.}
According to (\ref{minimal-power-function}) and  (\ref{power-function}),
$p^{\beta T}$ is a sum of Jacobi polynomials of degree $n$, then by using
 the orthogonality of the Jacobi polynomials
it is easy to obtain
\begin{equation} \label{int-zero}
\int_0^1 \tau^{l-1} \,p^{\beta T}(\tau)
\,d\tau=0, \ \text{ for any } l \in \{1,\dots,n\}.
\end{equation}
Then by applying $(\ref{int-zero})$, (\ref{mean2}) and (\ref{covariance2})  with the conditions given in $(\ref{condition1})$ and
$(\ref{condition2})$ we obtain
\begin{equation}
E\left[\ErrEstim{\beta T}{\varpi}(t_0)\right]=\int_0^1 p^{\beta T}(\tau)\,
E\left[\varpi(\beta T \tau)\right] \,d\tau.\label{mean}
\end{equation}
\begin{equation}
\begin{split}
&Cov\left[e_\varpi^{\beta T_1}(t_0), e_\varpi^{\beta T_2}(t_0)\right]\\=& \int_0^1 \int_0^1 p^{\beta T_1}(\tau)\,p^{\beta T_2}(s)\, Cov\left[ \varpi(\beta T_1 \tau),\, \varpi(\beta T_2
s)\right]\,d s \,d\tau. \label{covariance} \end{split}
\end{equation}

Consequently the mean and covariance  of $\ErrEstim{\beta T}{\varpi}(t_0)$ do
not depend on $t_0$. If we take $T_1=T_2$ in $(\ref{covariance})$,
then  the variance  of $\ErrEstim{\beta T}{\varpi}(t_0)$ do not depend on
$t_0$. Moreover, if $E[\varpi(\tau)]=\displaystyle\sum_{i=0}^{n-1}
\bar{\nu}_i \,\tau^i$, then by applying  $(\ref{int-zero})$ to
$(\ref{mean})$, we obtain $E[\ErrEstim{\beta T}{\varpi}(t_0)]=0$.
If $Cov[\varpi(s),\varpi(\tau)]=\displaystyle\sum_{i=0}^{n_1}
\bar{\eta}_i \,\tau^i\,\displaystyle\sum_{i=0}^{n_2} \bar{\eta}'_i \,
s^i$ with  $\min(n_1,n_2)\leq n-1$ then  by applying  $(\ref{int-zero})$ to
$(\ref{covariance})$, we obtain
$Cov[e_\varpi^{\beta T_1}(t_0),e_\varpi^{\beta T_2}(t_0)]=0$. Then if we take $T_1=T_2$ in $(\ref{covariance})$,
we get $Var[\ErrEstim{\beta T}{\varpi}(t_0)]=0$.
\hfill
$\blacksquare$

From which the following important theorem is obtained.

\begin{theorem}\label{theorem_casen}
Let $\ErrEstim{\beta T}{\varpi}(t_0)$ be the noise error contribution for
the Jacobi estimator $\Estim{\mu,\kappa}{\beta T,N,\xi}{x}^{(n)}({t_0})$ where the noise  $\{\varpi(\tau), \tau \geq 0\}$ satisfies
conditions ($C_1$) to ($C_3$), then
\begin{equation}
\ErrEstim{\beta T}{\varpi}(t_0)=0  \text { almost surely}.
\end{equation}
\end{theorem}

\noindent\textbf{Proof.} If the noise $\{\varpi(\tau), \tau \geq 0\}$ satisfies conditions ($C_1$) to ($C_3$), then we have
$E[\ErrEstim{\beta T}{\varpi}(t_0)]=0$ and $Var[\ErrEstim{\beta T}{\varpi}(t_0)]=0$. Since
$$E\left[\left(\ErrEstim{\beta T}{\varpi}(t_0)\right)^2\right]=Var\left[\ErrEstim{\beta T}{\varpi}(t_0)\right]+\left(E
\left [\ErrEstim{\beta T}{\varpi}(t_0)\right]\right)^2,$$
we get $E\left[\left(\ErrEstim{\beta T}{\varpi}(t_0)\right)^2\right]=0$.
Consequently, we have  $\ErrEstim{\beta T}{\varpi}(t_0)=0$ almost surely.

\hfill $\blacksquare$

Two stochastic processes, the Wiener process (also known as the
Brownian motion) and the Poisson process (cf \cite{SP}), play a
central role in the theory of stochastic processes. These processes
are valuable, not only as models of many important phenomena, but
also as building blocks to model other complex  stochastic
processes. They are characterized by:
\begin{itemize}
  \item \noindent let $\{W(t), t\geq 0\}$ be the Wiener process with
parameter $\sigma^2$, then
\begin{equation} \label{Wiener}
E\left[W(t)\right]=0, \ Cov\left[W(t),W(s)\right]= \sigma^2
\min(t,s);
\end{equation}
  \item let $\{N(t), t\geq 0\}$ be the Poisson process with intensity $\nu
\in \mathbb{R}^+$, then
\begin{equation}\label{poisson}
E\left[N(t)\right]=\nu t, \ Cov\left[N(t),N(s)\right]= \nu
\min(t,s).
\end{equation}
\end{itemize}
Thus, these processes satisfy  conditions $(C_1)$ and $(C_2)$. Hence,  we can characterize  the noise error
contributions due to these two stochastic processes for the Jacobi
estimators, and calculate the corresponding means and variances. If
the noise is a Wiener process, then it is clear that
$E[\ErrEstim{\beta T}{\varpi}(t_0)]=0$. If the noise is a
Poisson process, then we have
\begin{proposition}
The means of the noise error contributions due to a Poisson process
for the Jacobi estimators are given by
\begin{equation}
\begin{split}
\left\{
  \begin{array}{ll}
    E\left[e_{\varpi,\beta T, N, \xi}^{\mu,\kappa}(t_0)\right]=0, & \text{ if } n \geq 2, \\
   E\left[\ErrEstim{\mu,\kappa}{\varpi,\beta T}(t_0)\right]=E\left[\ErrEstim{\mu,\kappa}{\varpi,\beta T,N=2,\xi}(t_0)\right]=\nu, & \text{ if }
   n=1,
  \end{array},
\right.
\end{split}
\end{equation}
where $\ErrEstim{\mu,\kappa}{\varpi,\beta T}(t_0)$ (resp.
$\ErrEstim{\mu,\kappa}{\varpi,\beta T,N=2,\xi}(t_0)$ ) is the noise
error contribution for the minimal Jacobi estimators
$\Estim{\mu,\kappa}{\beta T}\dot x(t_0)$ (resp.  the
affine Jacobi estimators $\Estim{\mu,\kappa}{\beta T,N,\xi}\dot x(t_0)$) which is the estimates of the first order
derivative of $x$.
\end{proposition}

\noindent\textbf{Proof.} For $n \geq 2$, this can be simply proved
by using Theorem \ref{corollary}. Thus we only need to compute the
means of the noise error contributions for the estimates of
$\dot{x}$. Let $n=1$ in $(\ref{alien})$, then the minimal
estimators can be written in the following form
\begin{equation} \label{1.31}
\begin{split}
\Estim{\mu,\kappa}{\beta T}{\dot{x}}(t_{0})
=\int_{0}^{1} p^{\mu,\kappa}_{\beta T}(\tau)\ y( \beta T\tau + t_0)\ d\tau,\\
\end{split}
\end{equation}
where $p^{\mu,\kappa}_{\beta T}(\tau)=\frac{1}{\beta
T}\frac{\Gamma(\mu+\kappa+4)}{\Gamma(\kappa+2)
\Gamma(\mu+2)}\left((\mu+\kappa+2)\tau-(\kappa+1)\right)(1-\tau)^{\mu}
\tau^{\kappa}.$

The affine estimators $\Estim{\mu,\kappa}{\beta T, 2,\xi}{\dot{x}}(t_{0})$ are given by  $(\ref{16})$
\begin{equation*}%
\Estim{\mu,\kappa}{\beta T, 2,\xi}{\dot{x}}(t_{0}) = \lambda
_1(\xi,\kappa,\mu)\Estim{\mu+1,\kappa}{\beta T}{\dot{x}}(t_{0})+\lambda_0(\xi,\kappa,\mu)\Estim{\mu,\kappa+1}{\beta T}{\dot{x}}(t_{0}),
\end{equation*}
where $\lambda_1(\xi,\kappa,\mu)=(\kappa+3)-(\mu+\kappa+5)\xi$ and $\lambda_0(\xi,\kappa,\mu)=1-\lambda_1(\xi,\kappa,\mu)$
were obtained in \cite{num0}. According
to (\ref{1.31}), it reads
\begin{equation} \label{1.32}
\Estim{\mu,\kappa}{\beta T,2,\xi}\dot{x}(t_0) =\int_0^1
p^{\mu,\kappa}_{\beta T,2,\xi}(\tau)\, y(\beta T\tau+ t_0) \, d\tau,
\end{equation}
where $p^{\mu,\kappa}_{\beta T,2,\xi}(\tau) =\lambda_1(\xi,\kappa,\mu)
p^{\mu+1,\kappa}_{\beta T}(\tau)+\lambda_0(\xi,\kappa,\mu)
p^{\mu,\kappa+1}_{\beta T}(\tau)$.

According to $(\ref{mean})$ we obtain
\begin{equation*}
E\left[\ErrEstim{\mu,\kappa}{\varpi,\beta T}(t_0)\right]= \nu \beta
T\int_{0}^{1} \tau\,  p^{\mu,\kappa}_{\beta T}(\tau) \ d\tau.
\end{equation*}
By using integration by parts and the classical Beta function, we
obtain
\begin{equation*}
E\left[\ErrEstim{\mu,\kappa}{\varpi,\beta T}(t_0)\right]= \nu.
\end{equation*}
Moreover, since $\lambda_0(\xi,\kappa,\mu)+\lambda_1(\xi,\kappa,\mu)=1$ , one gets
\begin{eqnarray*}
E\left[\ErrEstim{\mu,\kappa}{\varpi,\beta T,2,\xi}(t_0)\right]&=&\lambda_1(\xi,\kappa,\mu)E\left[\ErrEstim{\mu+1,\kappa}{\beta T}(t_0)\right]+\lambda_0(\xi,\kappa,\mu)E\left[\ErrEstim{\mu,\kappa+1}{\beta T}(t_0)\right]\\
&= &\nu.
\end{eqnarray*}
Thus, this proof is completed.
\hfill $\blacksquare$

Now, in order to get the error bounds for the noise error contributions using the Bienaym\'{e}-Chebyshev
$(\ref{B-Chebyshev})$ we should compute the variance.  Since the covariance kernels of the Wiener process and the Poisson
process are determined by the same function $\min(\cdot,\cdot)$, the variances of the noise error
contributions due to a Wiener process or a Poisson process for the
Jacobi estimators $\hat{D}^{\mu,\kappa}_{\beta T} x^{(n)}(t_0)$ is given by (Using $(\ref{covariance})$ with $T=T_1=T_2$)
\begin{equation*}
Var\left[\ErrEstim{\mu,\kappa}{\varpi,\beta T}(t_0)\right]
=\eta\int_0^1 \int_0^1 p^{\mu,\kappa}_{\beta T}(\tau)\,p^{\mu,\kappa}_{\beta T}(s)\, \min (\beta Ts,\, \beta T
\tau) \,ds \,d\tau,
\end{equation*}
where $p^{\mu,\kappa}_{\beta T}(\tau)=\frac{(-1)^n w_{\mu+n,\kappa+n}^{(n)}(\tau)}{(\beta T)^n B(\mu+n+1,\kappa+n+1)}$. Using the symmetry property of  function $\min(\cdot,\cdot)$ and the fact that $\int_\tau^1 p^{\mu,\kappa}_{\beta T}(s) \,d s
=-\int_0^\tau p^{\mu,\kappa}_{\beta T}(s) \,d s$, we
obtain
\begin{equation}\label{eq_var_gene}
Var\left[e^{\mu,\kappa}_{\varpi,\beta T}(t_0)\right] = 2\eta
T\int_0^1 p^{\mu,\kappa}_{\beta T}(\tau)\,\tau \int_\tau^1
p^{\mu,\kappa}_{\beta T}(s) \,d s \,d\tau.
\end{equation}
Since
\begin{equation}\label{eq_simple_int}
\begin{split}
&\int_0^1 w_{\mu+n,\kappa+n}^{(n)}(\tau) \tau \int_{\tau}^1 w_{\mu+n,\kappa+n}^{(n)}(s) ds\, d\tau
=\\
&n! (n-1)! \int_0^1 w_{2\mu+1,2\kappa+2}(\tau) \,P^{\mu,\kappa}_{n}(\tau)\, P^{\mu+1,\kappa+1}_{n-1}(\tau)\, d\tau
\end{split}
\end{equation}
we have
\begin{equation}\label{eq_var_simple}
Var\left[e^{\mu,\kappa}_{\varpi,\beta T}(t_0)\right] =
\frac{2\eta n! (n-1)!}{T^{2n-1}B^2(\kappa+n+1,\mu+n+1)} I(\mu,\kappa,n).
\end{equation}
with
\begin{equation}\label{eq_var_gene}
I(\mu,\kappa,n) = \int_0^1 w_{2\mu+1,2\kappa+2}(\tau) \,P^{\mu,\kappa}_{n}(\tau)\, P^{\mu+1,\kappa+1}_{n-1}(\tau)\, d\tau.
\end{equation}
Let us stress that $Var\left[e^{\mu,\kappa}_{\varpi,\beta T}(t_0)\right]\sim \frac{1}{T^{2n-1}}$.

\par For $n=1$, we have the following results:
\begin{proposition}
The variances of the noise error contributions for the Jacobi
estimators of the first order derivative of ${x}$ are given by
\begin{equation} \label{variance-minimal}
Var\left[\ErrEstim{\mu,\kappa}{\varpi,\beta T}(t_0)\right]
=\frac{2\eta}{T}\frac{\mu+1}{2\mu+2\kappa+5}
\frac{B(2\mu+2,2\kappa+3)}{B^2(\kappa+2,\mu+2)},
\end{equation}
for minimal estimators and by
\begin{equation} \label{variance-affine}
\begin{split}
&Var\left[\ErrEstim{\mu,\kappa}{\varpi,\beta T,2,\xi}(t_0)\right]=\\
&\lambda_1^2(\xi,\kappa,\mu)\frac{2\eta}{T}\frac{\mu+2}{2\mu+2\kappa+7}
\frac{B(2\mu+4,2\kappa+3)}{B^2(\kappa+2,\mu+3)}\\&+\lambda_0^2(\xi,\kappa,\mu)\frac{2\eta}{T}\frac{\mu+1}{2\mu+2\kappa+7}
\frac{B(2\mu+2,2\kappa+5)}{B^2(\kappa+3,\mu+2)}\\
&+\lambda_0(\xi,\kappa,\mu) \lambda_1(\xi,\kappa,\mu)\frac{2\eta}{T}
\frac{B(2\mu+4,2\kappa+4)}{B(\kappa+2,\mu+3)B(\kappa+3,\mu+2)}
\end{split}
\end{equation}
for affine estimators. The value  $\eta$ is  equal to $\sigma^2$, if
the noise is a Wiener process, and $\eta$ is equal to $\nu$, if the
noise is a Poisson process.
\end{proposition}

\noindent\textbf{Proof.} Since
\begin{equation*}
I(\mu,\kappa,n=1)
=\frac{(\mu+1) B(2\mu+2,2\kappa+3)}{2\mu+2\kappa+5},
\end{equation*}
and using (\ref{eq_var_simple}) one gets the desired result. Similarly for
$Var\left[\ErrEstim{\mu,\kappa}{\varpi,\beta T,2, \xi}(t_0)\right]$.
\hfill $\blacksquare$\\

As a consequence, since $E[\ErrEstim{\mu,\kappa}{\varpi,\beta T}(t_0)]=E[\ErrEstim{\mu,\kappa}{\varpi,\beta T2,\xi}(t_0)]=\nu$ for a Wiener process ($\nu=0$) or Poisson process ($\nu\neq0$, where $\nu$ is the intensity parameter of the Poison Process), using the well known Bienaym\'{e}-Chebyshev
$(\ref{B-Chebyshev})$ we obtain the error bounds for the noise error contributions for the Jacobi estimators of the first order derivative of ${x}$.

\begin{theorem}[First order derivative estimation]\label{B-Chebyshev2} Let $n=1$. Let the noise be a Wiener process or Poisson process, then for any real number
$\gamma>0$,
\begin{equation} \label{B-Chebyshev_bis}
Pr\left(\left|\ErrEstim{\mu,\kappa}{\varpi,\beta T}(t_0)]-\nu\right| < \gamma\sqrt{Var[\ErrEstim{\mu,\kappa}{\varpi,\beta T}(t_0)]}\right) > 1-
\frac{1}{\gamma^2},
\end{equation}
\begin{equation} \label{B-Chebyshev_ter}
Pr\left(\left|\ErrEstim{\mu,\kappa}{\varpi,\beta T,2,\xi}(t_0)]-\nu\right| < \gamma\sqrt{Var[\ErrEstim{\mu,\kappa}{\varpi,\beta T,2,\xi}(t_0)]}\right) > 1-
\frac{1}{\gamma^2},
\end{equation} where $\nu=0$ for a Wiener process; $\nu\neq0$ for a Poisson process and
$Var[\ErrEstim{\mu,\kappa}{\varpi,\beta T}(t_0)]$,
$Var[\ErrEstim{\mu,\kappa}{\varpi,\beta T,2,\xi}(t_0)]$ are given
respectively by  $(\ref{variance-minimal})$ and
$(\ref{variance-affine})$.
\end{theorem}

\begin{figure}[h!]
 \centering
 \subfigure[$n=1$]
 {\includegraphics[scale=0.39]{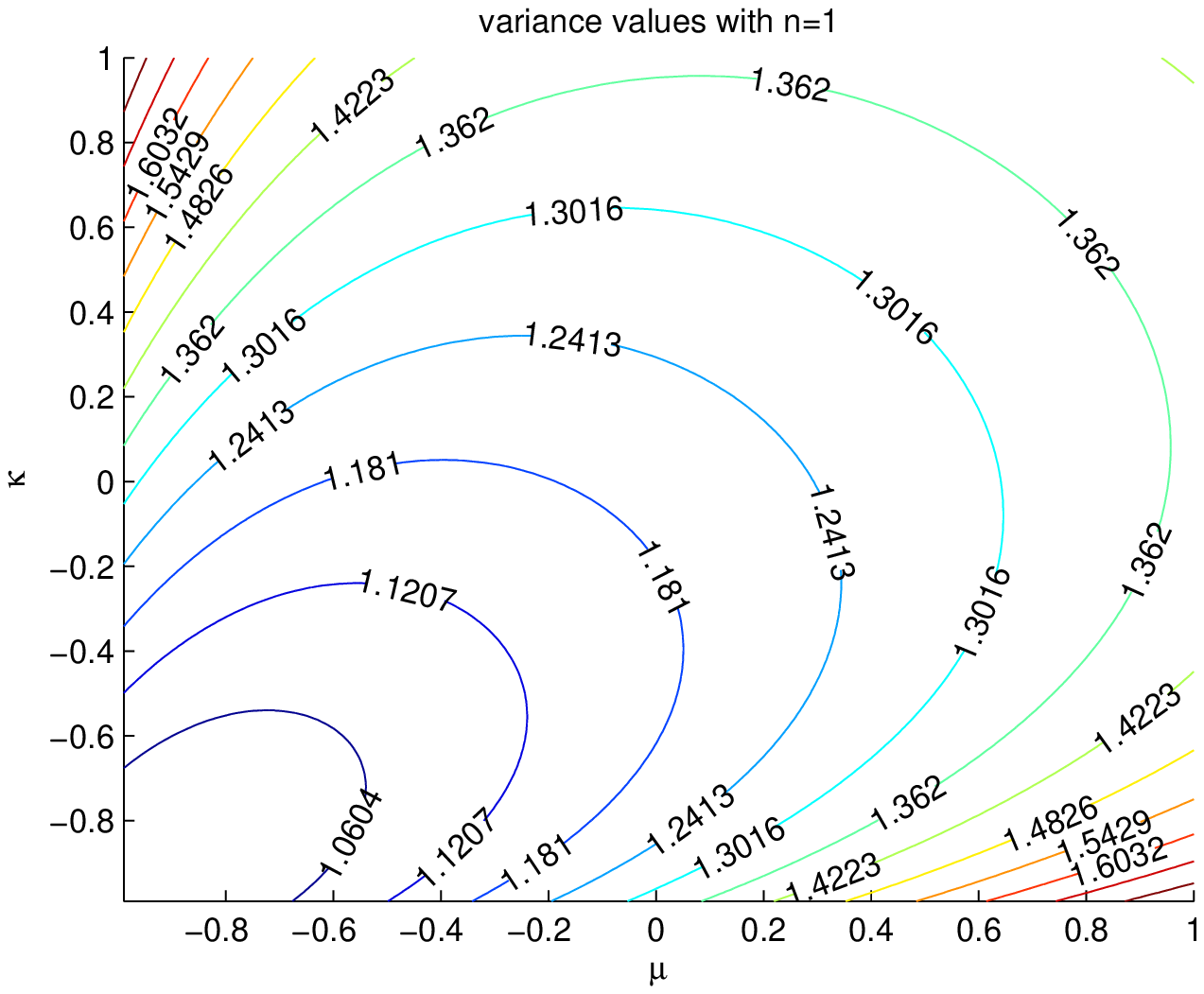}}
  \subfigure[$n=2$]
 {\includegraphics[scale=0.39]{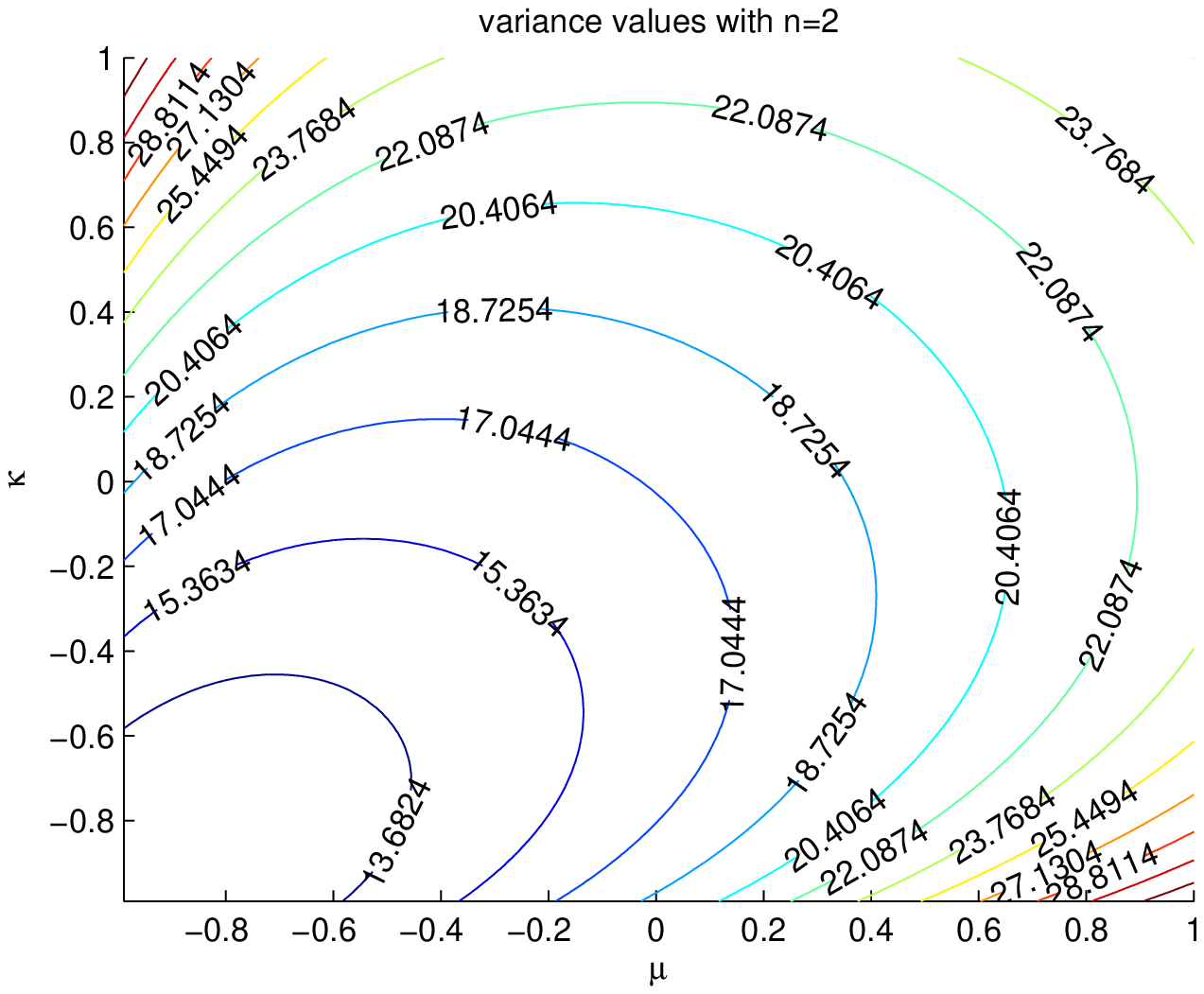}}
 \subfigure[$n=3$]
 {\includegraphics[scale=0.39]{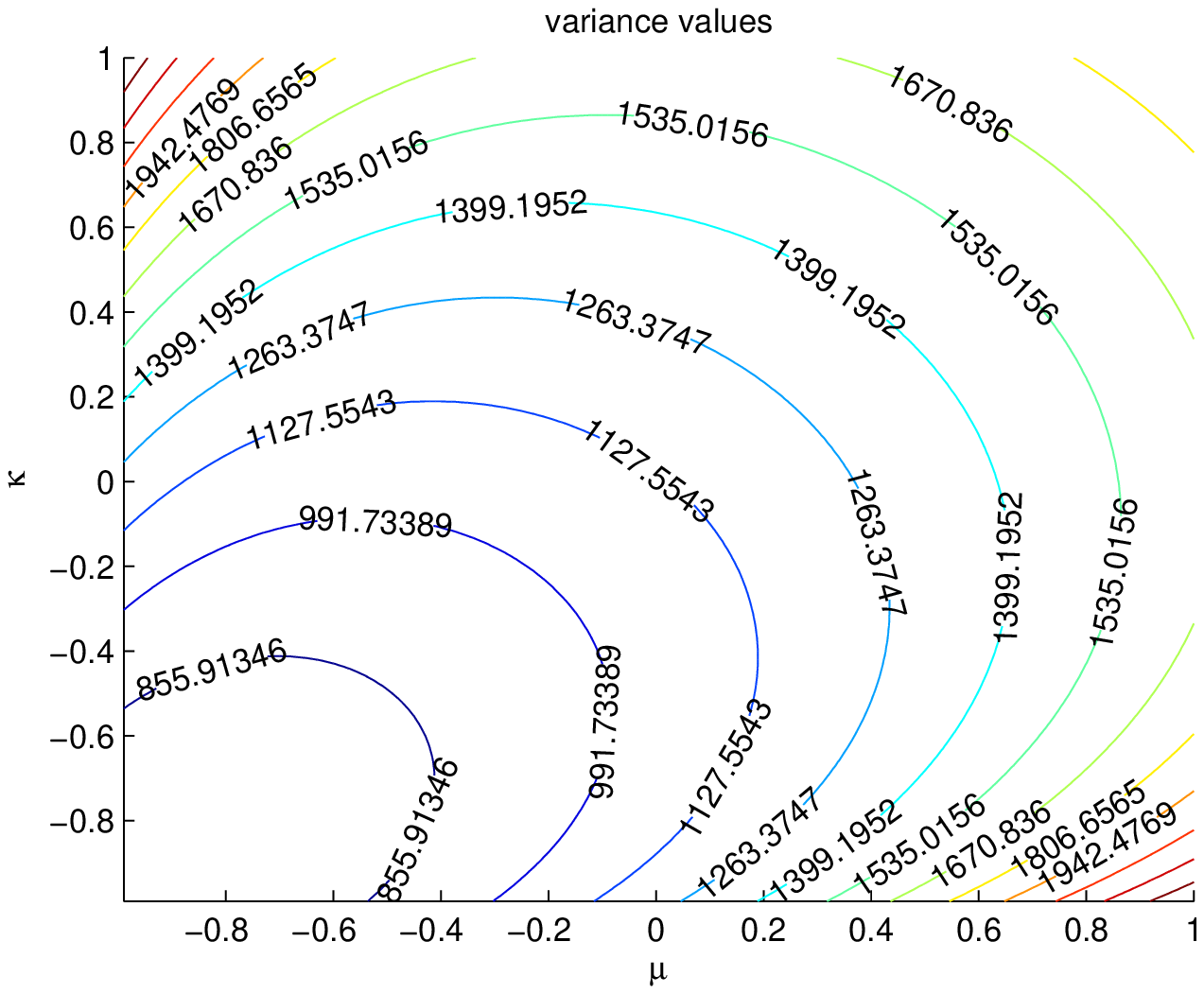}}
  \subfigure[$n=4$]
 {\includegraphics[scale=0.39]{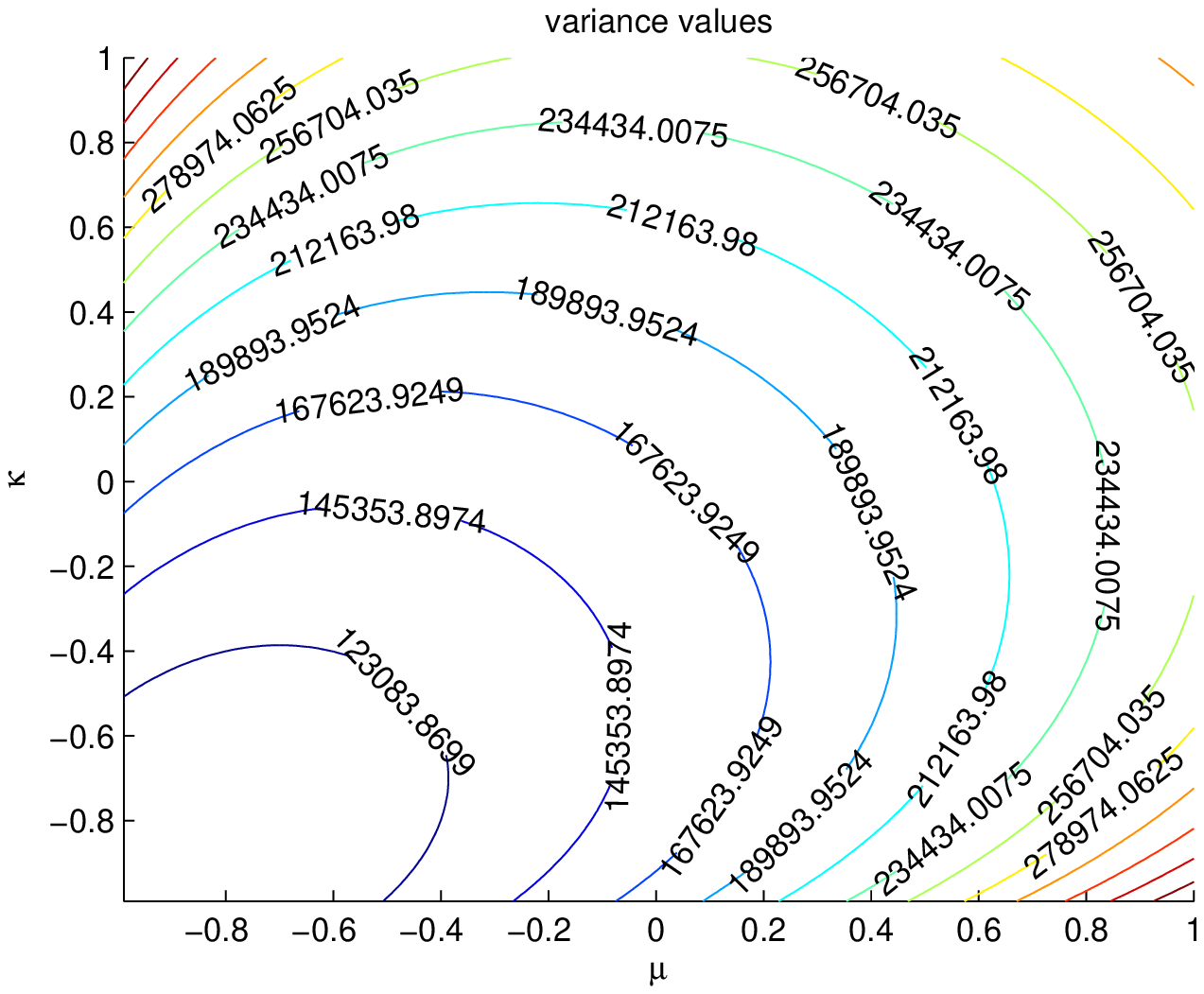}}
\caption{Variances of the noise errors for the minimal estimators.}
\label{figure_variance_min}
\end{figure}

\begin{figure}[h!]
 \centering
 \subfigure[$n=1$]
 {\includegraphics[scale=0.39]{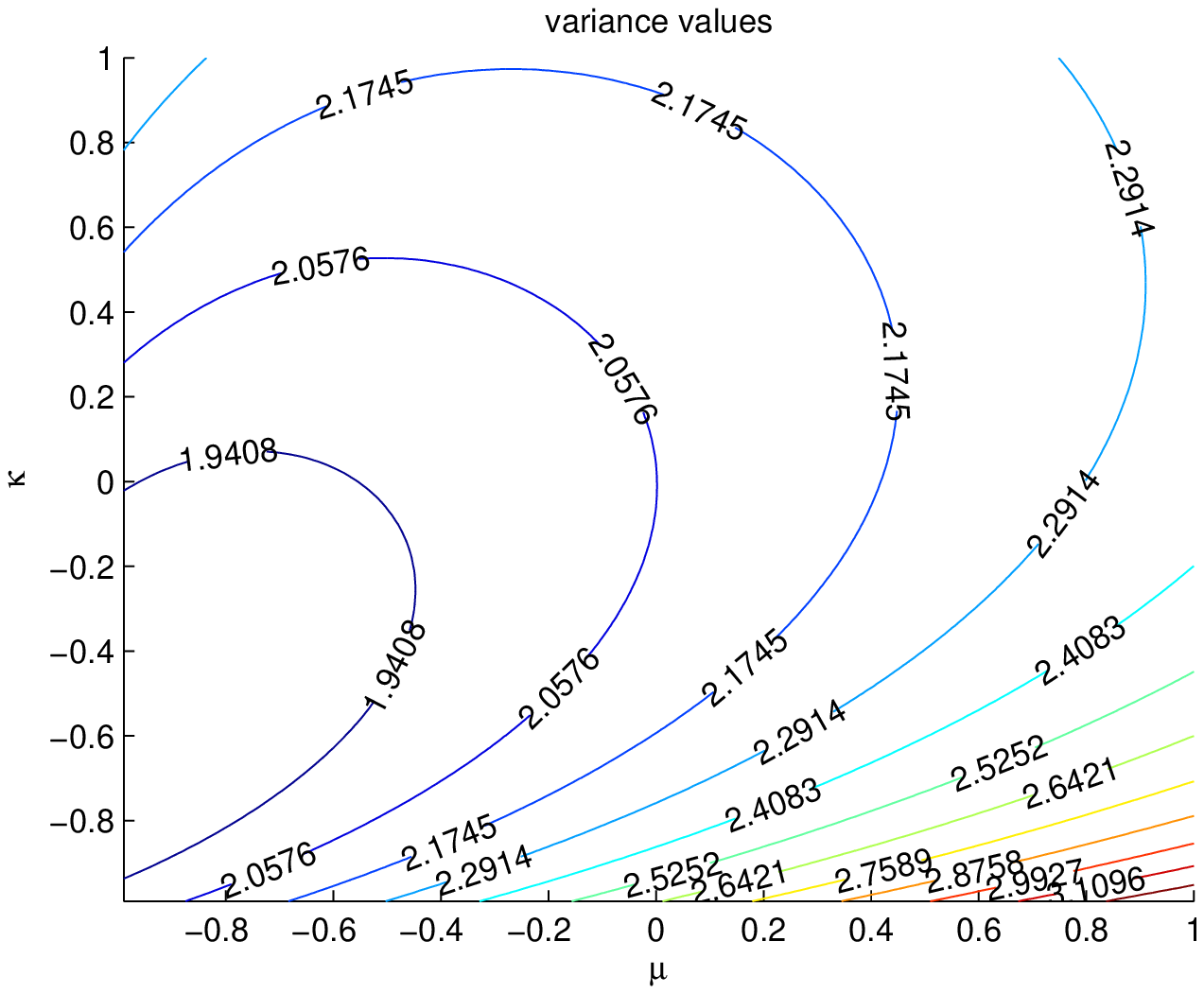}}
  \subfigure[$n=2$]
 {\includegraphics[scale=0.39]{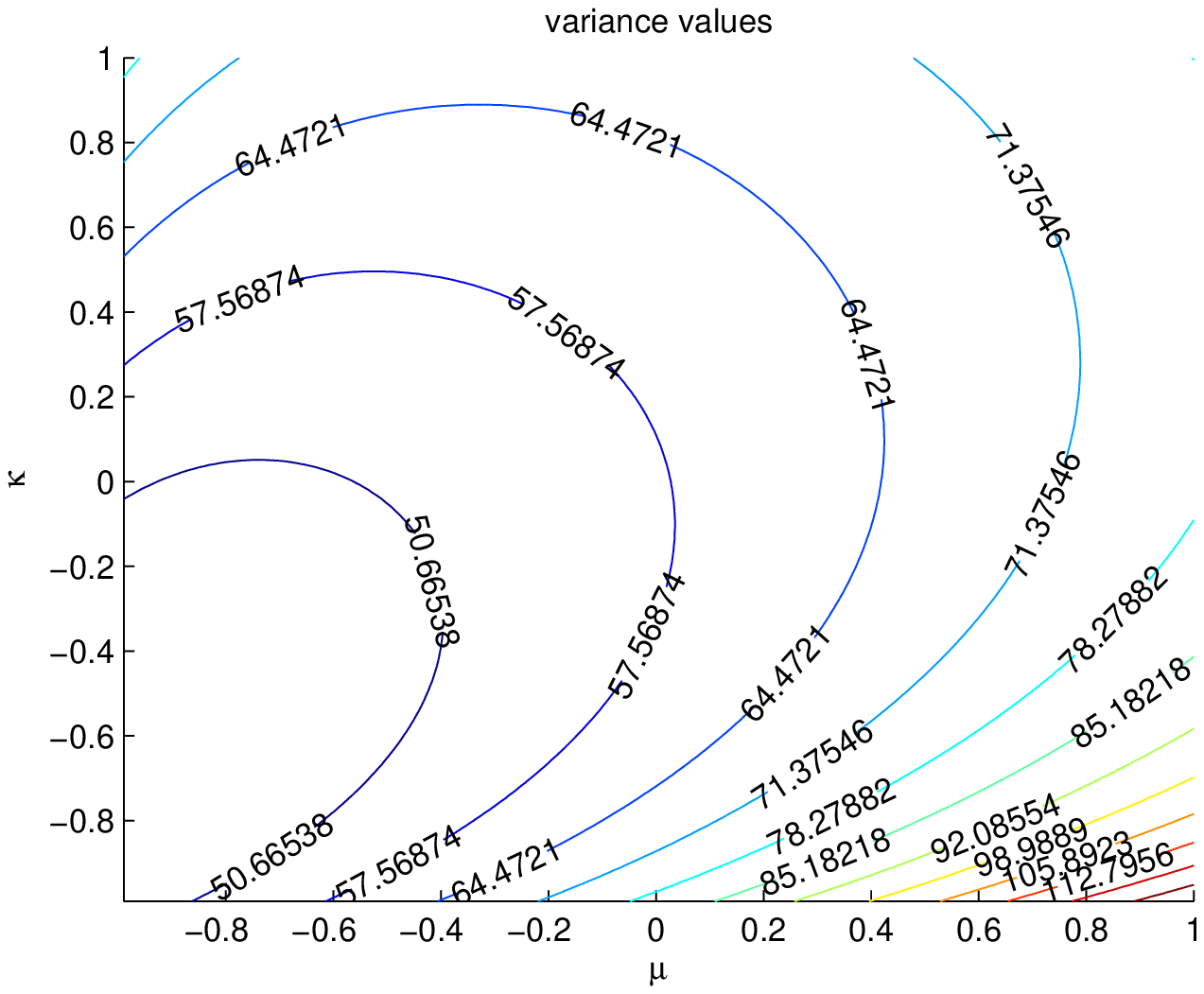}}
 \subfigure[$n=3$]
 {\includegraphics[scale=0.39]{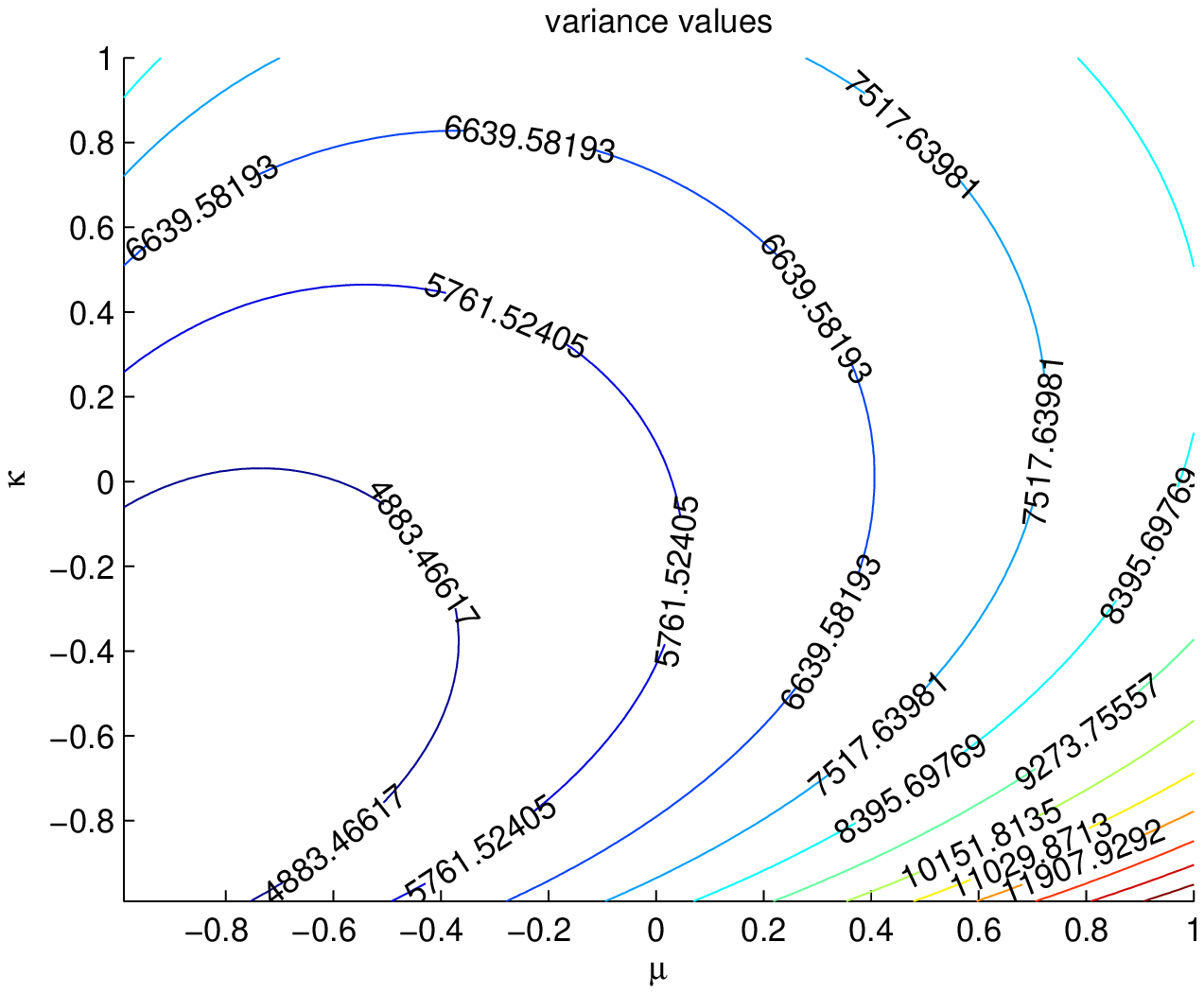}}
  \subfigure[$n=4$]
 {\includegraphics[scale=0.39]{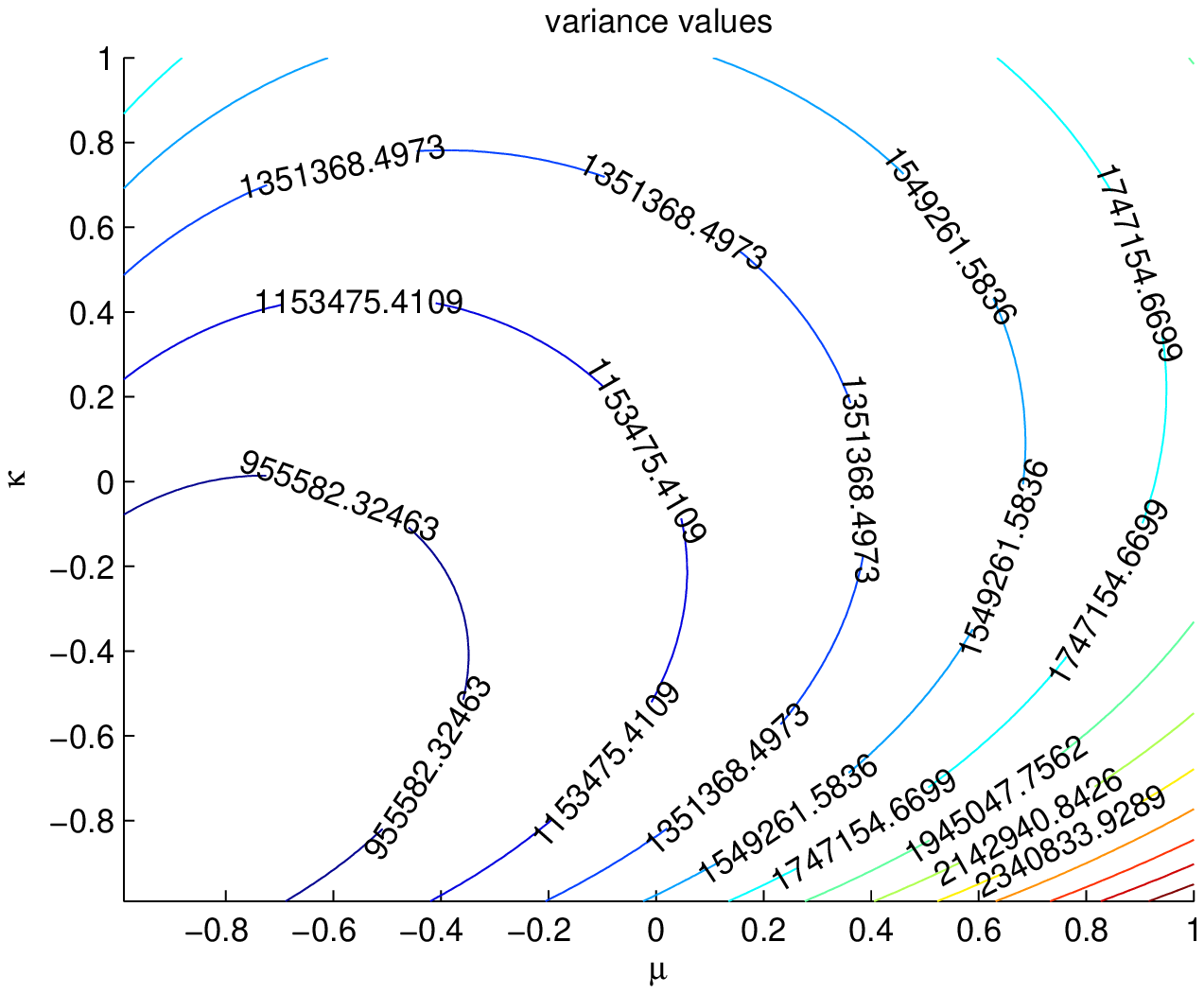}}
\caption{Variances of the noise errors for the affine estimators.}
\label{figure_variance_aff}
\end{figure}

%
For the case $n=1$, the bounds given by Theorem $\ref{B-Chebyshev2}$
characterize the noise error contribution
$\ErrEstim{\mu,\kappa}{\varpi,-T}(t_0)$ (respectively
$\ErrEstim{\mu,\kappa}{\varpi,-T,2,\xi}(t_0)$) for the Jacobi
estimation $\Estim{\mu,\kappa}{\varpi,-T}\dot{x}(t_{0})$
(respectively $\Estim{\mu,\kappa}{\varpi,-T,2,\xi}\dot{x}(t_{0})$).
They depend on
$Var\left[\ErrEstim{\mu,\kappa}{\varpi,-T}(t_0)\right]$ given by
(\ref{variance-minimal}) (respectively
$Var\left[\ErrEstim{\mu,\kappa}{\varpi,-T,2,\xi}(t_0)\right]$ given
by (\ref{variance-affine})).
Similar results can be obtained for $n=2$ since
\begin{equation*}
\begin{split}
2I(\mu,\kappa,n=2)=&-{(\kappa+2)^2(\kappa+1)} B(2\mu+5,2\kappa+3)\\
&+{(\kappa+2)(\mu+2)(3\kappa+5)}B(2\mu+4,2\kappa+4) \\&- {(\kappa+2)(\mu+2)(3\mu+5)} B(2\mu+3,2\kappa+5) \\
&  +{(\mu+2)^2(\mu+1)} B(2\mu+2,2\kappa+6).
\end{split}
\end{equation*}
and of course for higher values of $n$. Remember  that, for fixed
$T$, we have $Var\left[e^{\mu,\kappa}_{\varpi,\beta
T}(t_0)\right]\sim \frac{1}{T^{2n-1}}$). Since all these variance
functions decrease with respect to $T$ independently of $\kappa$ and
$\mu$, it is sufficient to observe the influence of $\kappa$ and
$\mu$. In the minimal Jacobi estimator case one can get a direct
computation (result is reported in Figure
$\ref{figure_variance_min}$ by taking $\eta=T=1$) whereas in the
affine case it is not difficult to obtain a 3-D plot as in Figure
$\ref{figure_variance_aff}$ and  where $\eta=T=1$,
$\xi=\xi(\kappa,\mu)$ is the smaller root of
$P_{2}^{\mu+1,\kappa+1}$.  From this analysis, we should take
negative values for $\kappa$ and $\mu$ so as to minimize the noise
error contribution. Moreover, we can observe that the variance of
$\ErrEstim{\mu,\kappa}{\varpi,-1,2,\xi}\dot{x}(t_{0})$ is larger
than the one of $\ErrEstim{\mu,\kappa}{\varpi,-1}\dot{x}(t_{0})$ if
we take same value for $\kappa$ and $\mu$, hence we should take  the
value of $T$ for affine estimator
$\Estim{\mu,\kappa}{\varpi,-T,2,\xi}\dot{x}(t_{0})$ larger than the
one for $\Estim{\mu,\kappa}{\varpi,-T}\dot{x}(t_{0})$ so as to
obtain the same noise effect.

\clearpage

\par Usually, the observation function $y$ is only known on discrete
values. Consequently, in the next subsection we will study discrete
noise cases.

\subsection{Noise in the discrete case}\label{subsection_discrete}

Let us now assume  that $y(t_i)= x(t_i) + \varpi(t_i)$ is a noisy
measurement of $x$ in discrete case with an equidistant sampling
period  $T_s$, where  noise $\varpi$ is assumed to satisfy condition
$(C_1)$ given in Subsection $\ref{subsection_stochastic}$. Let us
recall that the Jacobi estimators of the $n^{th}$ derivative of $x$
can be rewritten as follows
\begin{equation}
\Estim{\mu,\kappa}{\beta T,N,\xi}x^{(n)}(t_0)=\int_0^1
p^{\mu\kappa}_{\beta T,n,N,\xi}(\tau)\, y(t_0 + \beta T\tau) \,
d\tau. \label{eq-estimator_generale2}
\end{equation}

Since $y$ is a discrete measurement, we need to use a
numerical integration method  to  approximate the integral value in
$(\ref{eq-estimator_generale2})$. Let $t_i= \frac{i}{m}$ and $w_i>0$
for $i=0,\dots, m$ with $m=\frac{T}{T_s} \in \mathbb{N}$ (except for
$w_0 \geq 0$ and $w_m \geq 0$)  be respectively the abscissas and
the weights for a given numerical integration method used in
$(\ref{eq-estimator_generale2})$. Weight $w_0$ (resp. $w_m$) is set
to zero in order to avoid the infinite values when $\kappa$ (resp.
$\mu$) is negative. Then, we have
\begin{equation}
\Estim{\mu,\kappa}{\beta T,N,\xi}x^{(n)}(t_0)\approx
\sum_{i=0}^m \frac{w_i}{m} \, p^{\mu,\kappa}_{\beta T,n,N,\xi}(t_i)\,
y(t_0 +\beta Tt_i). \label{integration_numerique}
\end{equation}
Hence, the noise error contribution $\ErrEstim{\mu,\kappa}{\varpi,\beta T,N,\xi}(t_0)$ can be written in  discrete cases as
follows
\begin{equation}
\ErrEstim{\mu,\kappa}{\varpi,\beta T,N,\xi}(t_0)= \sum_{i=0}^m
\frac{w_i}{m} \, p^{\mu,\kappa}_{\beta T,n,N,\xi}(t_i)\, \varpi(t_0
+\beta Tt_i). \label{noise-contribution-discret}
\end{equation}

This numerical integration method also implies  an error which will
be studied in a future work. Consequently  the Jacobi estimators
lead to
\begin{equation}\label{general333}
\Estim{\mu,\kappa}{\beta T,N,\xi}x^{(n)}(t_0)=x^{(n)}(t_0)+ e_{m}(t_0)+ \ErrEstim{\mu,\kappa}{R_n,\beta T,N,\xi,m}(t_0)+\ErrEstim{\mu,\kappa}{\varpi,\beta T,N,\xi,m}(t_0),
\end{equation}
where   $\ErrEstim{\mu,\kappa}{R_n,\beta T,N,\xi,m}(t_0)$ is the
bias term error in  discrete cases,
$\ErrEstim{\mu,\kappa}{\varpi,\beta T,N,\xi,m}(t_0)$ is  the noise
error contribution in discrete cases (which will be shortly denoted
by $\ErrEstim{\beta T}{\varpi,m}(t_0)$ hereafter) and $e_{m}(t_0)$
is the numerical integration  error. To simplify the notations, as
in the previous section, $p^{\beta T}$ denotes power function
$p^{\mu,\kappa}_{\beta T, n,N,\xi}$. Then by applying the properties
of the mean, variance  and  covariance, we have
\begin{equation}
E\left[\ErrEstim{\beta T}{\varpi,m}(t_0)\right]= \frac{1}{m}\sum_{i=0}^m w_i
\, p^{\beta T}(t_i)\, E\left[\varpi(t_0 +\beta Tt_i)\right],
\label{mean-discret}
\end{equation}
\begin{equation} \label{variance-discret}
\begin{split}
&Var\left[\ErrEstim{\beta T}{\varpi,m}(t_0)\right]=
\frac{1}{m^2}\sum_{i=0}^m w_i^2 \, (p^{\beta T}(t_i))^2\,
Var\left[\varpi(t_0 +\beta
Tt_i)\right]\\&+\frac{2}{m^2}\sum_{i=0}^{m-1} \sum_{j=i+1}^{m}
w_i\,w_j\, p^{\beta T}(t_i)\,p^{\beta T}(t_j)\, Cov\left[\varpi(t_0
+\beta Tt_i),\varpi(t_0 +\beta Tt_j)\right].
\end{split}
\end{equation}
Moreover, for any $T_1>0$ and $T_2>0$
\begin{equation} \label{covariance-discret}
\begin{split}
&Cov\left[\ErrEstim{\beta T_1}{\varpi,m}(t_0),\ErrEstim{\beta T_2}{\varpi,m}(t_0)
\right]\\=&\frac{1}{m^2}\sum_{i=0}^{m} \sum_{j=0}^{m} w_i\,w_j\,
p^{\beta T_1}(t_i)\,p^{\beta T_2}(t_j)\, Cov\left[\varpi(t_0 +\beta
T_1 t_i),\varpi(t_0 +\beta T_2 t_j)\right].
\end{split}
\end{equation}

Now, by using Bienaym\'{e}-Chebyshev $(\ref{B-Chebyshev})$ and the
previous formulae, we can derive similar results than the ones
obtained in the previous subsection and which coincide if
$m\rightarrow \infty$. However this is true with some few additional
assumptions as detailed below.

In order to show the bridge with the previous Subsection
\ref{subsection_stochastic}, we will use the following properties,
where $T$, $T_1$ and $T_2$ are given (finite), and  $T_s$ tends to
$0$, i.e. $m$ tends to infinite.
\begin{equation}
\lim_{m \rightarrow \infty} E\left[\ErrEstim{\beta T}{\varpi,m}(t_0)\right]
=E\left[\ErrEstim{\beta T}{\varpi}(t_0)\right]
,\label{limite1}
\end{equation}
\begin{equation}
\lim_{m \rightarrow \infty} Var\left[\ErrEstim{\beta T}{\varpi,m}(t_0)\right]=Var\left[\ErrEstim{\beta T}{\varpi}(t_0)\right]
,\label{limite2}
\end{equation}
\begin{equation} \label{limite4}
\lim_{m \rightarrow \infty} Cov\left[\ErrEstim{\beta T_1}{\varpi,m}(t_0),
\ErrEstim{\beta T_2}{\varpi,m}(t_0) \right]
=Cov\left[\ErrEstim{\beta T_1}{\varpi}(t_0),
\ErrEstim{\beta T_2}{\varpi}(t_0)\right].
\end{equation}

From now on, let us consider a family of noises which are continuous
parameter stochastic processes satisfying  the following conditions
\begin{description}
  \item[$(C_4):$] the mean value and variance functions of $\{\varpi(\tau), \tau \geq 0\}$ are
continuous functions;
  \item[$(C_5):$] for any $s,t \geq 0$, $s\neq t$, $\varpi(s)$ and $\varpi(t)$ are
  independent.
\end{description}
Note that  white Gaussian noise and Poisson noise satisfy these conditions. Then, we can give the following theorem.

\begin{theorem} \label{theorem} Let $\{\varpi(\tau), \tau \geq 0\}$ be a continuous parameter stochastic process satisfying conditions
$(C_4)$ and $(C_5)$. Let $\varpi(t_i)$ be a sequence of $\{\varpi(\tau), \tau \geq 0\}$
with an equidistant sampling period $T_s$.
 If $\kappa, \mu>-\frac{1}{2}$, then we have
\begin{equation}\label{limite3}
\lim_{m \rightarrow \infty} Var\left[\ErrEstim{\beta
T}{\varpi,m}(t_0)\right]=0,
\end{equation}
where $\ErrEstim{\beta
T}{\varpi,m}(t_0)$ is the associated noise error
contribution for the Jacobi estimators defined by
(\ref{noise-contribution-discret}).
\end{theorem}

\noindent\textbf{Proof.} Since $\varpi(t_i)$ is a sequence of
independent random variables, by using $(\ref{variance-discret})$ we
have

\begin{equation}
\begin{split}
Var\left[\ErrEstim{\beta
T}{\varpi,m}(t_0)\right]=&
\frac{1}{m^2}\sum_{i=0}^m w_i^2 \, (p^{\beta T}(t_i))^2\,
Var\left[\varpi(t_0 +\beta Tt_i)\right].
\end{split}
\end{equation}
Since the variance function of $\varpi$ is continuous, we have
\begin{equation} \label{inegalite}
0\leq\frac{1}{m^2}\sum_{i=0}^m w_i^2 \, (p^{\beta T}(t_i))^2\,
\left|Var\left[\varpi(t_0 +\beta Tt_i)\right]\right| \leq
U\frac{w(m)}{m}\sum_{i=0}^m \frac{w_i}{m} \, (p^{\beta T}(t_i))^2,
\end{equation}
where $w(m)=\displaystyle\max_{0\leq i \leq m} w_i$ and
$U=\displaystyle\sup_{0 \leq t \leq 1}\left|Var\left[\varpi(t_0
+\beta Tt)\right]\right| < \infty$. Moreover,
\begin{equation} \label{}
\lim_{m \rightarrow \infty}\sum_{i=0}^m \frac{w_i}{m} \, (p^{\beta T}(t_i))^2=\int_0^1 (p^{\beta T}(t))^2 \,d t.
\end{equation}
Since $p^{\mu,\kappa}_{\beta T, n, N, \xi}$ is a sum of
$p^{\mu,\kappa}_{\beta T, n}$ (see $(\ref{power-function})$),
according to the expression of $p^{\mu,\kappa}_{\beta T, n}$ given
in $(\ref{minimal-power-function})$, $\int_0^1 (p^{\beta T}(t))^2
\,d t < \infty$ if $\int_0^1 (1-\tau)^{2\mu} \tau^{2k} \, d\tau <
\infty$. Consequently, as all $w_i$ are bounded, if $\kappa,
\mu>-\frac{1}{2}$ then
$$\displaystyle\lim_{m \rightarrow \infty}
U\frac{w(m)}{m}\displaystyle\sum_{i=0}^m \frac{w_i}{m} \, (p^{\beta T}(t_i))^2= 0.$$
The proof  is completed. \hfill
$\blacksquare$

According to the previous theorem, if we apply the
Bienaym\'{e}-Chebyshev inequality and (\ref{limite1}), then we can
obtain that $\ErrEstim{\beta
T}{\varpi,m}(t_0)$ converges in probability to
$\int_0^1 p^{\beta T}(\tau)\, E[\varpi(t_0 + \beta T\tau)] \, d\tau$
when $T_s \rightarrow 0$. Moreover, if we use the fact that
$E\left[\left(Y_{m}-c\right)^2\right]=Var\left[Y_{m}\right]+\left(E
\left [Y_{m}\right]-c\right)^2$ for any sequence of random variables
$Y_m$, then we can get the convergence in mean square.

\begin{corollary}
Let $\{\varpi(\tau), \tau \geq 0\}$ be a continuous parameter stochastic process satisfying conditions $(C_4)$ and $(C_5)$.
Let $\varpi(t_i)$ be a sequence of $\{\varpi(\tau), \tau \geq 0\}$
with an equidistant sampling period $T_s$.
 If $\kappa, \mu>-\frac{1}{2}$,  then $\ErrEstim{\beta
T}{\varpi,m}(t_0)$ converges in mean square to $\int_0^1 p^{\beta T}(\tau)\,
E[\varpi(t_0 + \beta T\tau)] \, d\tau$
 when $T_s \rightarrow 0 $,
where $p^{\beta T}$ is defined in $(\ref{power-function})$.
Moreover, if $E[\varpi(\tau)]=\displaystyle\sum_{i=0}^{n-1}
\bar{\nu}_i \,\tau^i$ with $\bar{\nu}_i \in \mathbb{R}$, then
$\ErrEstim{\beta
T}{\varpi,m}(t_0)$ converges in mean square to $0$
 when $T_s \rightarrow 0$.
\end{corollary}

\noindent\textbf{Proof.} If
$E[\varpi(\tau)]=\displaystyle\sum_{i=0}^{n-1} \bar{\nu}_i \,\tau^i$
with $\bar{\nu}_i \in \mathbb{R}$, then similarly  to Theorem
\ref{corollary} we can obtain  $\int_0^1 p^{\beta T}(\tau)\,
E[\varpi(t_0 + \beta T\tau)] \, d\tau=0$. Hence, this proof is
completed.

\hfill $\blacksquare$

\section{Numerical experiments} \label{section04}

If $\kappa$ (resp. $\mu$) is negative, $p^{\mu,\kappa}_{\beta T, n,
N, \xi}$ may be infinite at $\tau=0$ (resp. $\tau=1$). In the
previous Subsection \ref{subsection_discrete}  we choose $w_0=0$
(resp. $w_m=0$) so as to avoid this problem. But this choice of
$w_0$ and $w_m$ implies an error (see \cite{lyness}). In order to
reduce this error, we replace the term $\tau^{\kappa}$ at $\tau=0$
(resp. $(1-\tau)^{\mu}$ at $\tau=1$) in $p^{\mu,\kappa}_{\beta T, n,
N, \xi}$ by ${\left(\frac{F}{m}\right)}^\kappa$  (resp.
${\left(\frac{F}{m}\right)}^\mu$)
 with {$F \in ]0,1]$}  when $\kappa$ (resp. $\mu$) is negative.
For example, the power  function in the minimal Jacobi estimators $\Estim{\mu,\kappa}{\beta T}\dot{x}(t_{0})$ is $$p^{\mu,\kappa}_{\beta T}(\tau)=\frac{1}{\beta
T}\frac{\Gamma(\mu+\kappa+4)}{\Gamma(\kappa+2)
\Gamma(\mu+2)}\left((\mu+\kappa+2)\tau-(\kappa+1)\right)(1-\tau)^{\mu}
\tau^{\kappa}.$$
If $\kappa<0$ and $\mu \geq 0$, then we take
\begin{equation}
p_{\kappa,\mu,\beta T}(0) \approx \frac{1}{\beta
T}\frac{\Gamma(\mu+\kappa+4)}{\Gamma(\kappa+2)
\Gamma(\mu+2)}\left(-(\kappa+1)\right)
\left({\frac{F}{m}}\right)^{\kappa}.
\end{equation}

In the two following subsections, we use the trapezoidal
rule as the numerical integration method.

\subsection{Simulation results with a Brownian motion noise}
\begin{figure}[h!]
\centering \subfigure[The noisy measurement
$y(t_i)$ and $x(t_i)$. ]  {\includegraphics[scale=0.6]{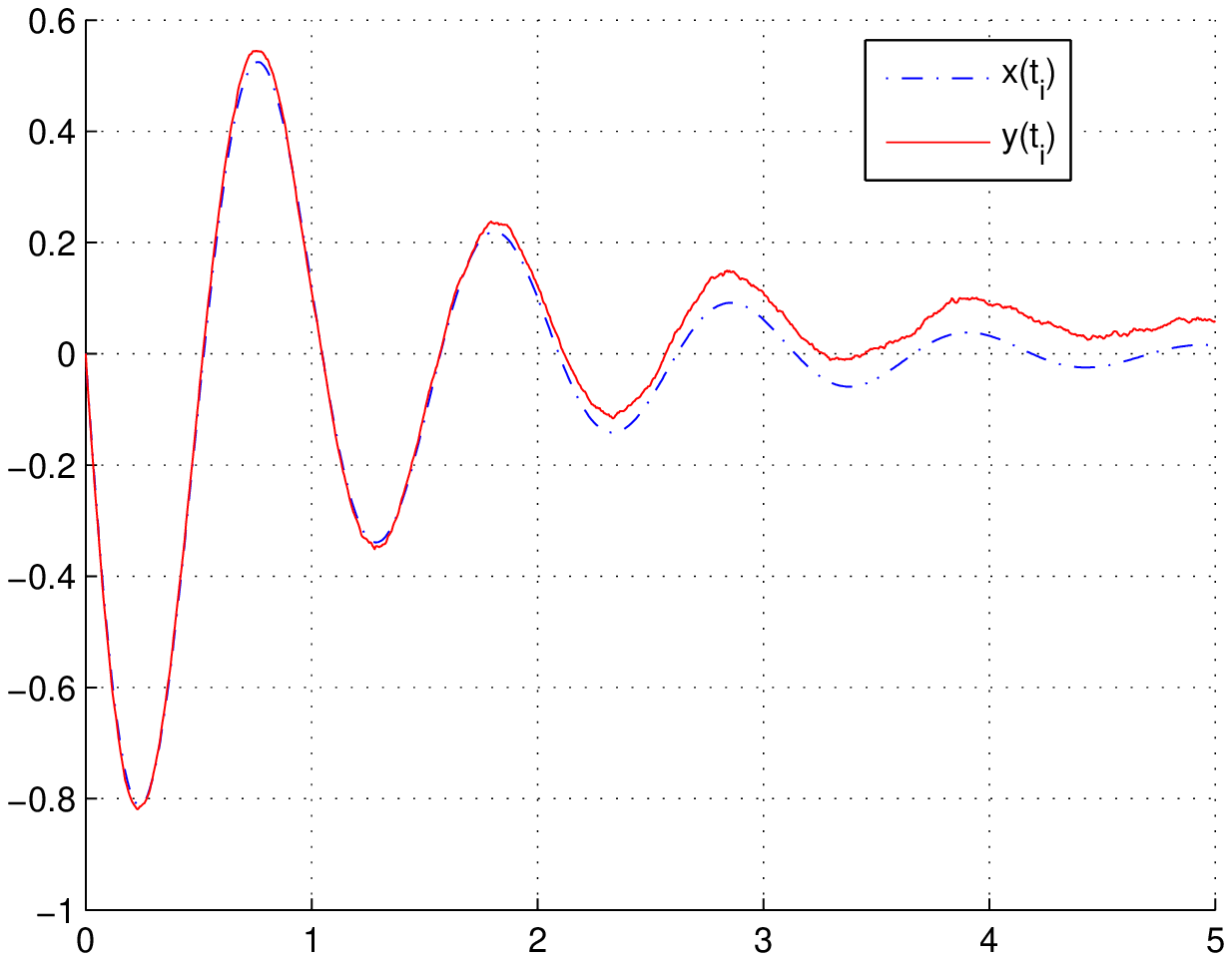}} %
\subfigure[ The associated Brownian noise $C \varpi(t_i)$.] %
{\includegraphics[scale=0.6]{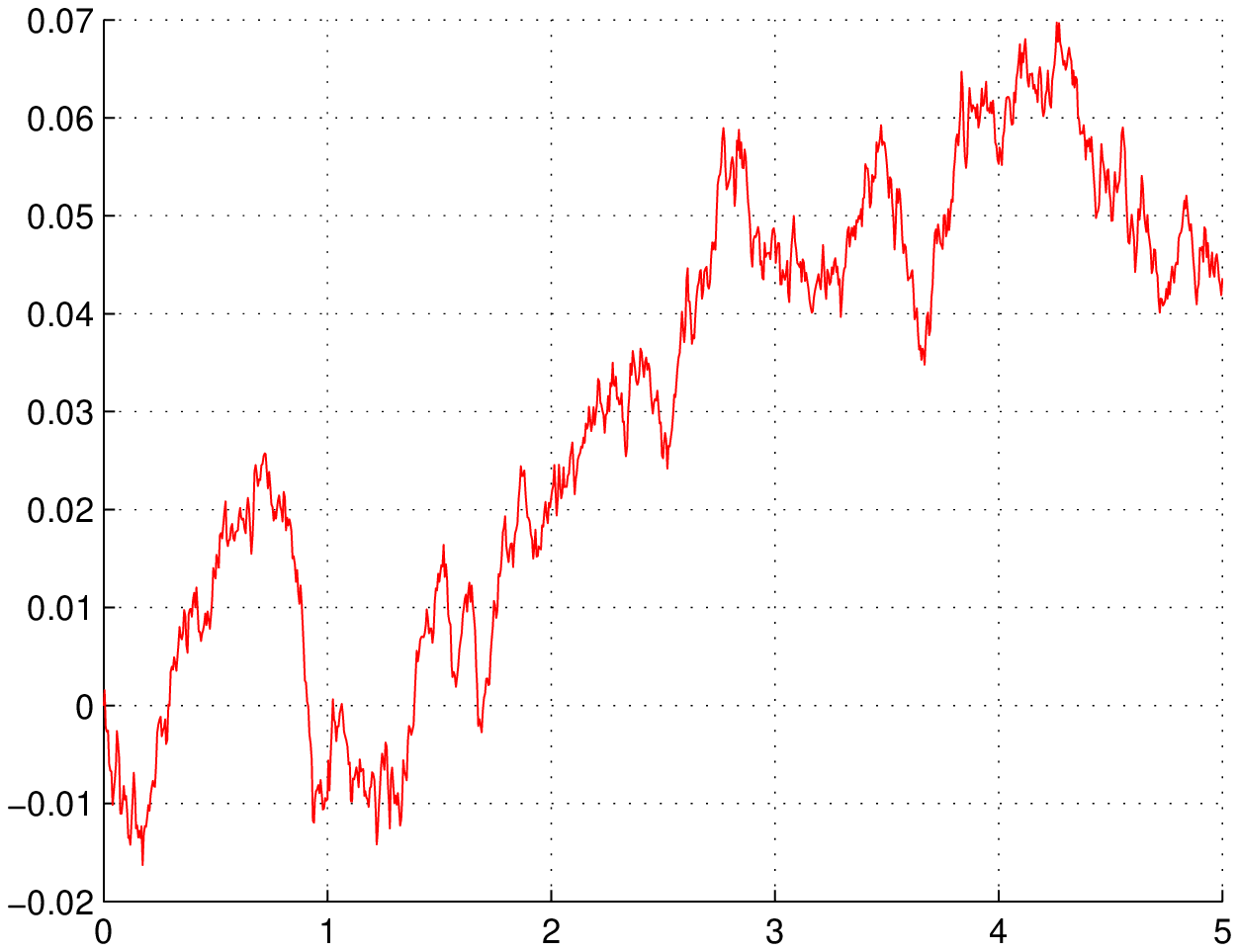}} \caption{The given data.} %
\label{figure_signal}
\end{figure}

\begin{figure}[h!]
 \centering
 \subfigure[ $\Estim{\mu,\kappa}{-T}\dot{x}(t_{0})$ with $\kappa=\mu=0$, $T=18T_s$ and
 $\Estim{\mu,\kappa}{-T}\dot{x}(t_{0})$ with $\kappa=-0.79$, $\mu=0$, $T=30T_s$ and $F=0.1$. ]
 {\includegraphics[scale=0.6]{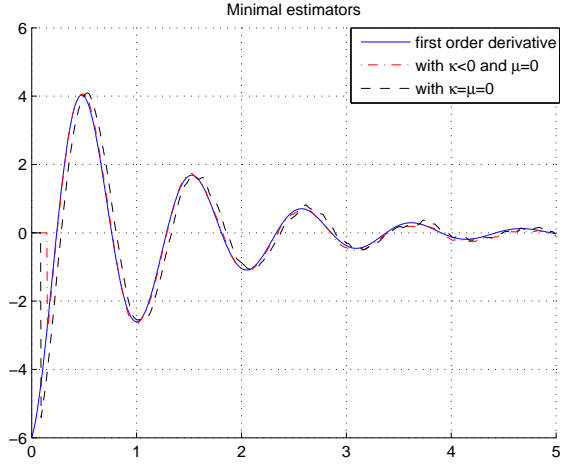}} \label{figure_mini}
 \subfigure[ $\Estim{\mu,\kappa}{-T,2, \xi}\dot{x}(t_{0})$ with $\kappa=\mu=0$, $T=30T_s$, $\xi=0.276$ and
 $\Estim{\mu,\kappa}{-T,2, \xi}\dot{x}(t_{0})$ with $\kappa=-0.78$, $\mu=-0.6$, $T=46T_s$, $\xi=0.218$ and $F=0.1$.]
 {\includegraphics[scale=0.6]{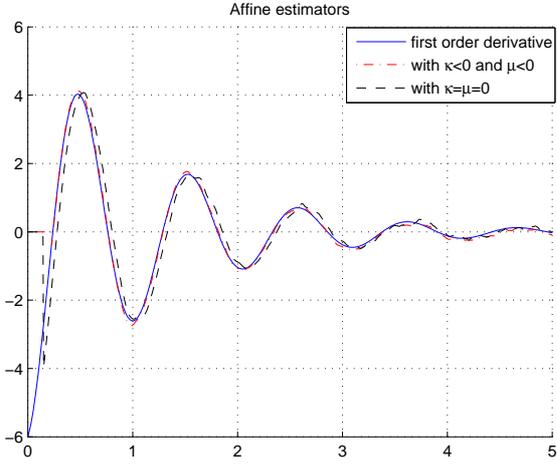}} \label{figure_affine}
\caption{Estimations obtained by using the Jacobi causal
estimators.} \label{figure_estimations}
\end{figure}

In this subsection, we assume  that $y(t_i)= x(t_i) + C \varpi(t_i)$ with $t_i=T_s
i$ for $i=0,\cdots,1000$ ($T_s=\frac{1}{200}$), is a noisy
measurement of  $x(t_i)= \exp(\frac{-t_i}{1.2}) \sin(6 t_i+\pi)$.
Noise $C \varpi$ is assumed to be a Brownian motion defined by
$(\ref{Wiener})$ with $\sigma^2=1$ and $C>0$. The coefficient $C$ is
chosen so that  the signal-to-noise ratio
$SNR=10\log_{10}\left(\frac{\sum |y(t_i)|^2}{\sum
|C\varpi(t_i)|^2}\right)$ is equal to $SNR=16 \text{dB}$ (see, e.g.,
\cite{R17} for this well known concept in signal processing). Figure $\ref{figure_signal}$ reports the noisy measurement with its associated noise.

We use the  minimal causal estimator
$\Estim{\mu,\kappa}{-T}\dot{x}(t_{0})$  defined in
$(\ref{1.31})$ and the affine causal  estimator
$\Estim{\mu,\kappa}{-T,2,\xi}\dot{x}(t_{0})$ defined in
 $(\ref{1.32})$ to estimate the first order derivative of ${x}$.
It was shown in the previous sections that  the bias term errors
(Section \ref{sec_biais_error}) and the noise error contributions
for the Jacobi estimators (Section \ref{section03}) both depend on
the parameters $\kappa, \mu, T$ and $\xi$. We have previously  shown
the parameters' influence on the time delay values (the bias term
errors) and the noise error contributions for the minimal estimators
$\Estim{\mu,\kappa}{-T}\dot{x}(t_{0})$ and the affine estimators
$\Estim{\mu,\kappa}{-T,2,\xi}\dot{x}(t_{0})$. In order to obtain an
``optimized'' error (minimum), it is clear that we should take a
negative value for $\kappa$. Concerning the choice of $T$ we should
make a compromise since small $T$ makes the
 bias term error small but also produces large noise error contribution. We take negative value
for $\mu$ so as  to reduce noise error contributions. We can see in
Figure $\ref{figure_estimations}$ the obtained estimations
respectively by using the minimal causal estimators  and
 the affine causal estimators.

In each figure, the solid line represents the exact derivative of
$x$, the dashed line represents the time-delayed estimation with
$\kappa, \mu \in \mathbb{N}$ and the dotted line represents the
estimation with $\kappa, \mu \in ]-1, +\infty[$. We can see  the
``delay-free''  estimations with  $\kappa, \mu \in ]-1, 0]$. In
fact, we can find out an appropriate value for $F$ and $m$, so that
the numerical integration method error reduces the bias term errors
and the time delay values. The analysis for such errors will be
studied in a future work.

\begin{figure}[h!]
 \centering
 \subfigure[Noise error contributions $\ErrEstim{\mu,\kappa}{\varpi,-T,2,\xi}(t_0)$ in  $\Estim{\mu,\kappa}{-T,2,\xi}\dot{x}(t_0)$
  with $\kappa=\mu=0$, $T=30T_s$, $\xi=0.276$. ]
 {\includegraphics[scale=0.6]{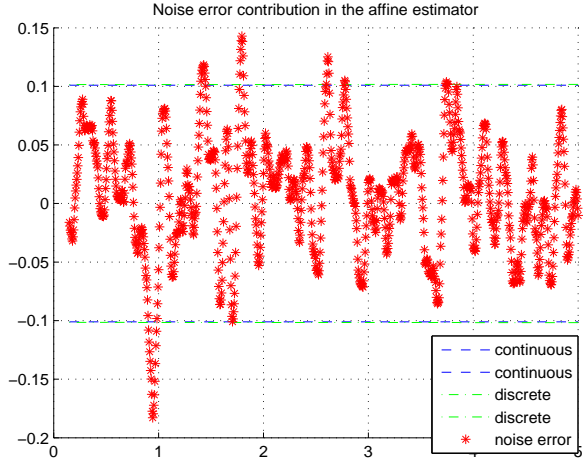}}  \label{errorbound1}
 \subfigure[ Noise error contributions $\ErrEstim{\mu,\kappa}{\varpi,-T,2,\xi}(t_0)$ in  $\Estim{\mu,\kappa}{-T,2,\xi}\dot{x}(t_0)$
 with $\kappa=-0.78$, $\mu=-0.6$, $T=46T_s$, $\xi=0.218$ and $F=0.1$.]
 {\includegraphics[scale=0.6]{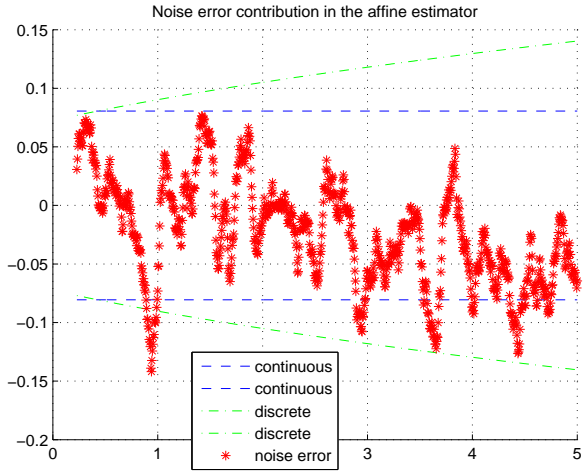}} \label{errorbound2}
\caption{Noise error contributions for the affine estimators and the
error bounds.}
 \label{errorbound}
\end{figure}

The associated noise error contributions for the estimations
obtained by using respectively the affine estimator with $\kappa,
\mu \in \mathbb{N}$ and the ones with $\kappa, \mu \in ]-1,+\infty[$
are shown in Figure $\ref{errorbound}$. In this figure, the error
bounds for these noise error contributions are also shown. The
dotted lines represent the error bounds obtained in the continuous
case by using Theorem
 $\ref{B-Chebyshev2}$ with $\gamma=2$. The dashed lines
represent the error bounds obtained in this discrete case by using a discrete counter part of Theorem
 $\ref{B-Chebyshev2}$ with $\gamma=2$.

 In order to
compare these estimations, we calculate the total error variance
$\int_0^5 e(\tau)^2 d\tau$ for each estimate. As the measurement is
discrete, we take $T_s\displaystyle\sum_{i=0}^{1000} e(\tau_i)^2$ as
the approximation of $\int_0^5 e(\tau)^2 d\tau$.  We also consider
the $SNR$ by calculating each estimate and the associated noise
error contributions. The error variance and the $SNR$ value for each
estimate are given in Table \ref{tabular2}. All the values are
calculated in the same time interval $[50T_s,5]$. We can see that
with the same $SNR$ value the estimations obtained with $\kappa, \mu
\in ]-1, 0]$ produce smaller total errors than the ones obtained
with $\kappa, \mu \in \mathbb{N}$. Moreover, we calculate the time
delay
 for the estimates obtained by using minimal causal estimator
with $\kappa, \mu \in \mathbb{N}$ (resp. affine causal estimator
with $\kappa, \mu \in \mathbb{N}$) which is given by
$\frac{\kappa+n+1}{\kappa+\mu+2n+2}T$ (resp. $T\xi(\kappa,\mu))$. We
will make the same comparison in the next subsection.

\begin{table}[h!]
\caption{Error variance, $SNR$ and time delay value in each
estimate} \label{tabular2}
\begin{tabular}{ccc}
\hline\noalign{\smallskip}
$F=0.1$ & $\Estim{0, 0}{-18T_s}\dot{x}(t_{0})$ & $ \Estim{0,-0.79}{-30T_s}\dot{x}(t_{0})$ \\
\noalign{\smallskip}\hline\noalign{\smallskip}
 $\int_0^5 e(\tau)^2 d\tau$ &$0.5506$ & $0.0181$  \\
$SNR$ & $27.8019$ &  $27.5795$  \\
$Theoretical\  Delay$ & $0.045$ & $0.0565$ \\
\hline\noalign{\smallskip} $F=0.1$ &$\Estim{0, 0}{-30T_s,2, 0.276}\dot{x}(t_{0})$  & $\Estim{-0.6,-0.78 }{-46T_s,2, 0.218}\dot{x}(t_{0})$\\
\hline\noalign{\smallskip}
$\int_0^5 e(\tau)^2 d\tau$   &  $0.5148$ & $0.0197$    \\
 $SNR$   & $27.6765$ & $27.9936$ \\
 $Theoretical\  Delay$   & $0.0414$ & $0.0501$ \\
 \hline\noalign{\smallskip}
\end{tabular}
\end{table}

\subsection{Simulations results with a white  Gaussian noise}
\begin{figure}[h!]
\centering \subfigure[The noisy measurement $y(t_i)$ and the smooth
function $x(t_i)$.] {\includegraphics[scale=0.6]{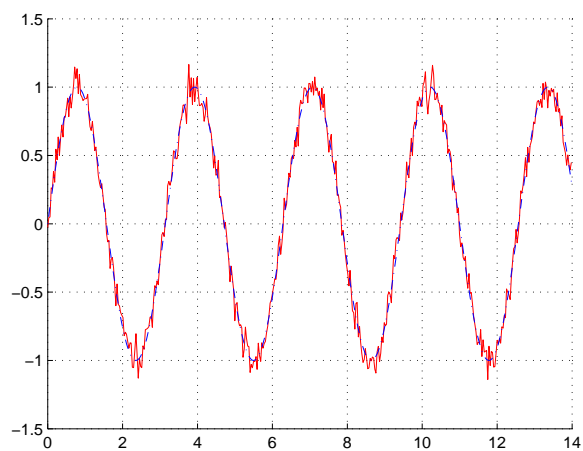}}
\subfigure[ The associated white Gaussian noise $C\varpi(t_i)$.]
{\includegraphics[scale=0.6]{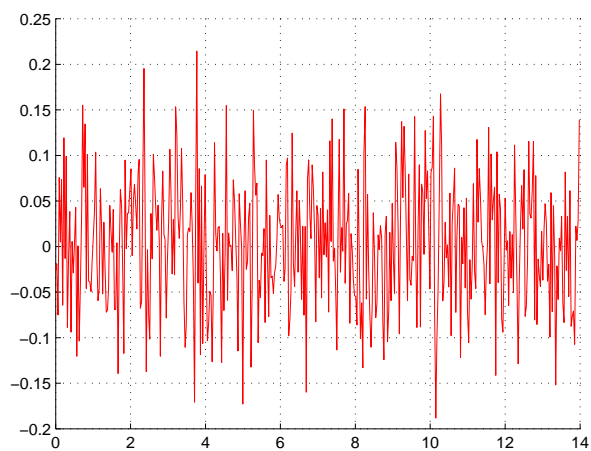}} \caption{The given data.}
\label{figure_signal2}
\end{figure}

Let $y(t_i)= \sin(2 t_i) + C\varpi(t_i)$, with $t_i=T_s i$ for
$i=0,\cdots,445$ ($T_s=\frac{\pi}{100}$), be a noisy measurement of
$\sin(2 t_i)$. The samples of noise $C \varpi(t_i)$ are simulated
from a zero-mean white Gaussian $iid$ sequence where coefficient $C$
is adjusted in such a way that $SNR=20 \text{dB}$ (see Figure
$\ref{figure_signal2}$).

\begin{figure}[h!]
 \centering
 \subfigure[ $\Estim{\mu,\kappa}{-T}\dot{x}(t_{0})$ with $\kappa=\mu=0$, $T=25T_s$ and
 $\Estim{\mu,\kappa}{-T}\dot{x}(t_{0})$ with $\kappa=-0.75$, $\mu=0$, $T=25T_s$ and $F=0.5$. ]
 {\includegraphics[scale=0.6]{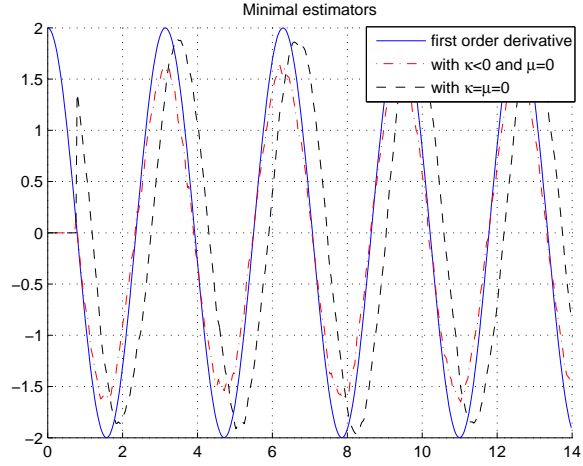}} \label{figure_mini1}
 \subfigure[ $\Estim{\mu,\kappa}{-T2, \xi}\dot{x}(t_{0})$ with $\kappa=\mu=0$, $T=38T_s$, $\xi=0.276$ and
 $\Estim{\mu,\kappa}{-T2, \xi}\dot{x}(t_{0})$ with $\kappa=-0.7$, $\mu=-0.66$, $T=32T_s$, $\xi=0.234$ and $F=0.5$.]
 {\includegraphics[scale=0.6]{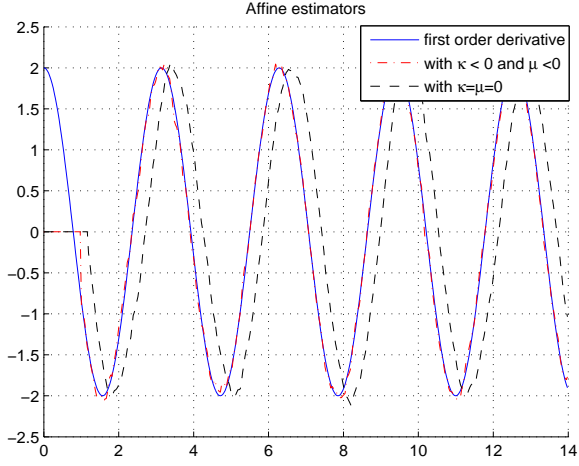}} \label{figure_affine1}
\caption{Estimations obtained by using the Jacobi causal
estimators.} \label{figure_estimations2}
\end{figure}

We use the  minimal causal estimator
$\Estim{\mu,\kappa}{-T}\dot{x}(t_{0})$  defined in
$(\ref{1.31})$ and the affine causal  estimator
$\Estim{\mu,\kappa}{-T,2,\xi}\dot{x}(t_{0})$ defined in
 $(\ref{1.32})$ to estimate the first order derivative of ${x}$.
We can see in Figure $\ref{figure_estimations2}$ the estimations
obtained respectively by using the minimal causal estimators  and
 the affine causal estimators.

In each figure,  the solid lines represent the exact derivative of
$x$, the dashed lines represent the time-delayed estimation with
$\kappa, \mu \in \mathbb{N}$ and the dotted lines represent the
``delay-free'' estimation with $\kappa, \mu \in ]-1, 0]$. The error
variance and the $SNR$ value for each estimate are given in Table
\ref{tabular1}. All the values are calculated in the same time
interval $[38T_s,14]$.

\begin{table}
\caption{Error variance, $SNR$ and time delay value in each
estimate} \label{tabular1}
\begin{tabular}{ccc}
\hline\noalign{\smallskip}
$F=0.5$ & $\Estim{0,0}{-25T_s}\dot{x}(t_{0})$ & $\Estim{0,-0.75}{-25T_s}\dot{x}(t_{0})$ \\
\noalign{\smallskip}\hline\noalign{\smallskip}
$\int_0^5 e(\tau)^2 d\tau$ &    $2.2351$ & $0.1855$  \\
$SNR$ & $27.5396$ & $27.7133$ \\
 $Theoretical\ Delay$ & $0.3927$ & $0.3021$\\
\hline\noalign{\smallskip}
$F=0.5$ & $\Estim{0,0}{-38T_s,2,0.276}\dot{x}(t_{0})$  & $\Estim{-0.66,-0.7}{-32T_s,2,0.234}\dot{x}(t_{0})$ \\
\noalign{\smallskip}\hline\noalign{\smallskip}
$\int_0^5 e(\tau)^2 d\tau$ &     $1.7919$ & $0.0085$ \\
$SNR$ & $27.5376$ & $27.2715$  \\
 $Theoretical\ Delay$ & $0.3295$ & $0.2352 $\\
\hline\noalign{\smallskip}
\end{tabular}
\end{table}

\section{Conclusion}

In this article, we study  recent algebraic parametric estimation
techniques introduced in \cite{num} which provide an estimate of the
derivatives by using iterated integrals of a noisy observation
signal. These algebraic parametric differentiation
techniques give derivative estimations which contain two sources of
errors: the bias term error and the noise error contribution. In
order to reduce these errors, we extend the parameter domains used
in the estimators. Then, we study some error bounds which depend on
these parameters. This allows us to minimize these errors. We show
that a compromise choice of these parameters  implies an
``optimized'' error among the noise error contribution, the bias
term error and the time delay. We also give some examples where the
errors due to numerical integration method permit us to further reduce the
time delay of the estimators. We will study this interesting fact in
a future work.



\end{document}